\documentclass[review]{classe}

\usepackage{lineno,hyperref}
\usepackage{placeins}
\usepackage{amsmath}
\usepackage{amsfonts}
\usepackage{amssymb}
\usepackage{gensymb}
\usepackage{subfigure}
\usepackage{color}
\usepackage{algorithm}
\usepackage{algorithmic}

\usepackage{pgfplots}
\pgfplotsset{compat=newest}
\pgfplotsset{plot coordinates/math parser=false}
\newlength\figureheight
\newlength\figurewidth
\pgfplotsset{
legend image code/.code={
\draw[mark repeat=2,mark phase=2]
plot coordinates {
(0cm,0cm)
(0.15cm,0cm)       
(0.32cm,0cm)         
};%
}
}
\usepackage{tikz}
\usetikzlibrary{arrows,calc,decorations.markings,decorations.pathreplacing,math,arrows.meta,fillbetween}
\tikzset{mark size=0.9} 
\usepackage{geometry}
\newcommand{\figref}[1]{Fig.~\ref{#1}}
\newcommand{\tabref}[1]{Tab.~\ref{#1}}
\newcommand{\secref}[1]{Sec.~\ref{#1}}

\newcommand{\algoref}[1]{Algorithm~\ref{#1}}

\newcommand{\st}{{\ \colon\ }}
\newcommand{\iii}{{\vert\kern-0.25ex\vert\kern-0.25ex\vert}}

\newcommand{\Div}{{\rm div}\,}

\newcommand{\eye}{{I}}
\newcommand{\velo}{\mathbf u}
\newcommand{\test}{\mathbf v}

\newcommand{\other}{\mathbf w}
\newcommand{\veloset}{{V}}
\newcommand{\pset}{{P}}
\newcommand{\ptest}{\pi}
\newcommand{\padj}{q}

\newcommand{\domvel}{\mathbf V}

\newcommand{\normal}{{\boldsymbol\nu}}
\newcommand{\normalcomp}{{\nu}}
\newcommand{\tangent}{{\boldsymbol\tau}}
\newcommand{\curv}{{H}}

\newcommand{\tn}{{(n)}}
\newcommand{\tnp}[1]{{(n+#1)}}
\newcommand{\tnm}[1]{{(n-#1)}}
\newcommand{\tN}{{(N)}}

\newcommand{\rd}[1]{%
  \IfEqCase{#1}{%
    {x}{\xi}%
    {y}{\eta}%
    {(x,y)}{(\xi,\eta)}%
  }[\PackageError{\rd}{Opzione non definita per rd: #1}{}]%
}

\newcommand{\zW}[2]{%
  \mathrel{\vbox{\offinterlineskip\ialign{%
    \hfil##\hfil\cr
    $\scriptscriptstyle\circ$\cr
    \noalign{\kern0.1ex}
    ${\,W^{#1}_{#2}}$\cr
}}}(\Omega^\gamma)}
\newcommand{\zWo}[2]{%
  \mathrel{\vbox{\offinterlineskip\ialign{%
    \hfil##\hfil\cr
    $\scriptscriptstyle\circ$\cr
    \noalign{\kern0.1ex}
    ${\,W^{#1}_{#2}}$\cr
}}}(\Omega^0)}
\newcommand{\zWsmall}[3]{%
  \mathrel{\vbox{\offinterlineskip\ialign{%
    \hfil##\hfil\cr
    $\scriptscriptstyle\circ$\cr
    \noalign{\kern0.1ex}
    $\scriptstyle{\,W^{#1}_{#2}}$\cr
}}}(\Omega^{#3})}

\newcommand{\Vo}{%
  \mathrel{\vbox{\offinterlineskip\ialign{%
    \hfil##\hfil\cr
    $\scriptscriptstyle\circ$\cr
    \noalign{\kern0.4ex}
    ${\,V_h}$\cr
}}}}

\newcommand{\veloadj}{\mathbf z}
\newcommand{\adjtest}{\mathbf w}
\newcommand{\padjtest}{\eta}
\newcommand{\sigmaadj}{{\varsigma}}

\newcommand{\control}{{\zeta}}
\newcommand{\controlvar}{{\boldsymbol\zeta}}
\newcommand{\controlset}{{\mathcal M_{ad}}}

\DeclareMathOperator*{\argmin}{{\rm arg\,min}}

\newcommand\mynote[1]{\textcolor{black}{#1}}

\newtheorem{thrm}{Theorem}[section]
\newtheorem{rmrk}[thrm]{Remark}

\begin{document}

\begin{frontmatter}

\title{Optimal control in ink-jet printing via
instantaneous control}

  \author{Ivan Fumagalli}

  \author{Nicola Parolini\corref{cor}}

  \author{Marco Verani}

\address{MOX, Dipartimento di Matematica, Politecnico di Milano, \\ P.za Leonardo da Vinci 32, I-20133 Milano, Italy}

\begin{abstract}
\mynote{This paper concerns the optimal control of a free surface flow with moving contact line, inspired by an application in ink-jet printing.
Surface tension, contact angle and wall friction are taken into account by means of the generalized Navier boundary condition.
The time-dependent differential system is discretized by an arbitrary Lagrangian-Eulerian finite element method, and a control problem is addressed by an instantaneous control approach, based on the time discretization of the flow equations.
The resulting control procedure is computationally highly efficient and its assessment by numerical tests show its effectiveness in deadening the natural oscillations that occur inside the nozzle and reducing significantly the duration of the transient preceding the attainment of the equilibrium configuration.}

\end{abstract}


\end{frontmatter}


\section{Introduction}\label{sec:intro}

Optimal control of systems governed by PDEs is a highly relevant issue in industrial applications, and the interest on this subject is continuously increasing in different fields, ranging from the mathematical description and analysis of this kind of problems to the engineering solutions to it, and to the actual implementation of control strategies in industrial processes.
This task is even more challenging if fluid dynamics systems are involved, due to the intrinsic complex nature of the related phenomena.

The present work is inspired by an application in inkjet printing.
This technology is widely employed to many aims, ranging from household usage to industry and security.
Particularly in the latter case, precision and accuracy of the printing are the main objectives, and thus the control of the ink jets plays an important role.

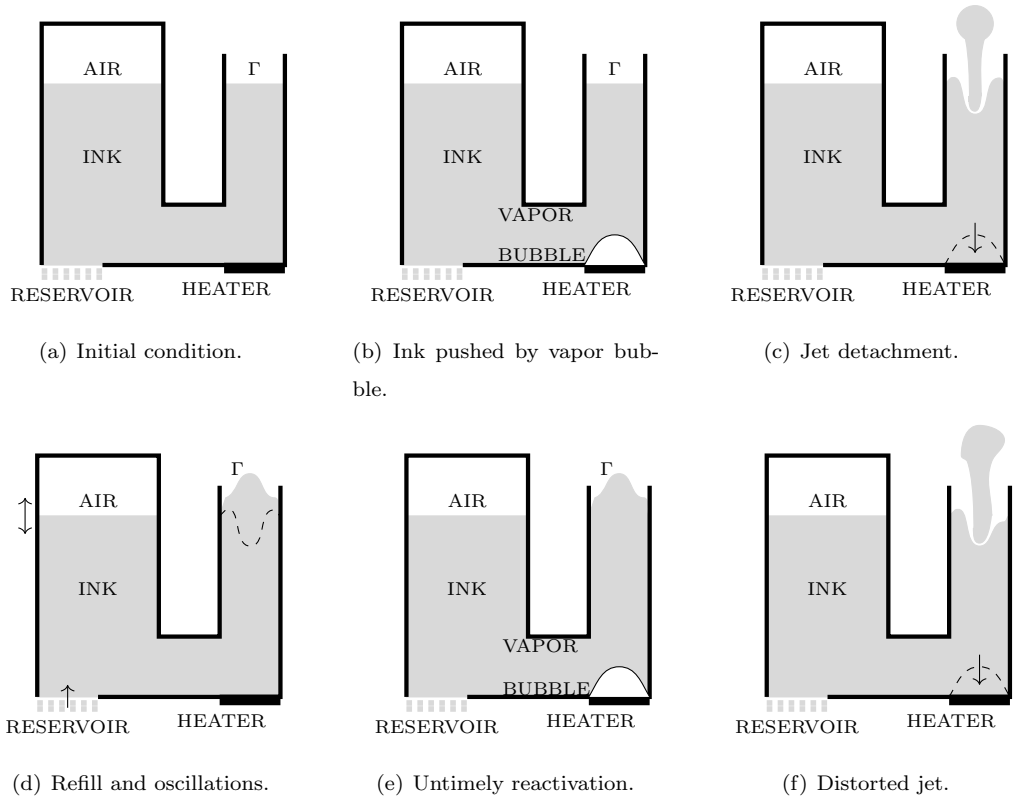
\begin{figure}[h]
\centering
\subfigure[Initial condition.
\label{fig:t0}]{
	\begin{tikzpicture}[scale=0.8]
		\draw [gray!30,fill=gray!30] (4,3)--(4,0)--(0,0)--(0,3)-- node[black,above]{\scriptsize AIR}(2,3)--(2,1)--(3,1)--(3,3)-- node[black,above]{\scriptsize $\Gamma$}(4,3);
		\draw [ultra thick] (4,3.5)--(4,0)--(1,0);
		\draw [ultra thick] (0,0)--(0,4)--(2,4)--(2,1)--(3,1)--(3,3.5);
		\node at (1,1.8) {\scriptsize INK};
		\draw [line width=4] (3,-0.05)-- node[black,at end, anchor=north east]{\scriptsize HEATER} (4,-0.05);
		\draw [gray!30, dotted, line width=2.1] (0,-0.1)--(1,-0.1);
		\draw [gray!30, dotted, line width=2.1] (0,-0.2)-- node[black,below]{\scriptsize RESERVOIR} (1,-0.2);
	\end{tikzpicture}
}\qquad%
\subfigure[Ink pushed by vapor bubble.
\label{fig:t1}]{
	\begin{tikzpicture}[scale=0.8]
		\draw [gray!30,fill=gray!30] (4,3)--(4,0)--(0,0)--(0,3)-- node[black,above]{\scriptsize AIR}(2,3)--(2,1)--(3,1)--(3,3)-- node[black,above]{\scriptsize $\Gamma$}(4,3);
		\draw [ultra thick] (4,3.5)--(4,0)--(1,0);
		\draw [ultra thick] (0,0)--(0,4)--(2,4)--(2,1)--(3,1)--(3,3.5);
		\node at (1,1.8) {\scriptsize INK};
		\draw [line width=4] (3,-0.05)-- node[black,at end, anchor=north east]{\scriptsize HEATER} (4,-0.05);
		\draw [gray!30, dotted, line width=2.1] (0,-0.1)--(1,-0.1);
		\draw [gray!30, dotted, line width=2.1] (0,-0.2)-- node[black,below]{\scriptsize RESERVOIR} (1,-0.2);
		\draw [fill=white] (3,0) to[out=60,in=180] (3.5,0.5) to     [out=0,in=120]
		    (4,0);
		\node[text width=1] at (1.6,0.5) {\scriptsize VAPOR\\[-1ex]\scriptsize BUBBLE};
	\end{tikzpicture}
}\qquad%
\subfigure[Jet detachment.
\label{fig:t2}]{
	\begin{tikzpicture}[scale=0.8]
		\draw [gray!30,fill=gray!30] (4,3)--(4,0)--(0,0)--(0,3)-- node[black,above]{\scriptsize AIR}(2,3)--(2,1)--(3,1)--(3,3) to (3.1,3.1) to[out=20,in=180] (3.5,2.5) to [out=0,in=180] (3.9,3.1) to [out=0,in=160] (4,3);
		\draw [gray!30,fill=gray!30] (3.5,4) circle (0.3cm);
		\draw [gray!30,fill=gray!30] (3.3,3.9) to [out=-70,in=-90] (3.4,2.8) to [out=-90,in=180] (3.5,2.55) to [out=0,in=90] (3.6,2.8) to [out=90,in=-110] (3.7,3.9)--(3.3,3.9);
		\draw [ultra thick] (4,3.5)--(4,0)--(1,0);
		\draw [ultra thick] (0,0)--(0,4)--(2,4)--(2,1)--(3,1)--(3,3.5);
		\node at (1,1.8) {\scriptsize INK};
		\draw [line width=4] (3,-0.05)-- node[black,at end, anchor=north east]{\scriptsize HEATER} (4,-0.05);
		\draw [gray!30, dotted, line width=2.1] (0,-0.1)--(1,-0.1);
		\draw [gray!30, dotted, line width=2.1] (0,-0.2)-- node[black,below]{\scriptsize RESERVOIR} (1,-0.2);
		\draw [dashed] (3,0) to[out=60,in=180] (3.5,0.5) to     [out=0,in=120]  (4,0);
		\draw [->] (3.5,0.7)--(3.5,0.2);
	\end{tikzpicture}
}\\
\subfigure[Refill and oscillations.
\label{fig:t3}]{
	\begin{tikzpicture}[scale=0.8]
		\draw [gray!30,fill=gray!30] (4,3)--(4,0)--(0,0)--(0,3)-- node[black,above]{\scriptsize AIR}(2,3)--(2,1)--(3,1)--(3,3) to (3.1,3.3) to[out=20,in=180] node[black,above]{\scriptsize $\Gamma$} (3.5,3.7) to [out=0,in=180] (3.9,3.3) to [out=0,in=160] (4,3);
		\draw [dashed] (3,3) to (3.1,3.1) to[out=20,in=180]     (3.5,2.5) to [out=0,in=180] (3.9,3.1) to [out=0,in=160] (4,3);
		\draw [ultra thick] (4,3.5)--(4,0)--(1,0);
		\draw [ultra thick] (0,0)--(0,4)--(2,4)--(2,1)--(3,1)--(3,3.5);
		\node at (1,1.8) {\scriptsize INK};
		\draw [line width=4] (3,-0.05)-- node[black,at end, anchor=north east]{\scriptsize HEATER} (4,-0.05);
		\draw [gray!30, dotted, line width=2.1] (0,-0.1)--(1,-0.1);
		\draw [gray!30, dotted, line width=2.1] (0,-0.2)-- node[black,below]{\scriptsize RESERVOIR} (1,-0.2);
		\draw [<->] (-0.2,2.7)--(-0.2,3.3);
		\draw [->] (0.5,-0.18)--(0.5,0.2);
	\end{tikzpicture}
	\begin{tikzpicture}[scale=0.8]
	\end{tikzpicture}
}\qquad%
\subfigure[Untimely reactivation.
\label{fig:t4}]{
	\begin{tikzpicture}[scale=0.8]
		\draw [gray!30,fill=gray!30] (4,3)--(4,0)--(0,0)--(0,3)-- node[black,above]{\scriptsize AIR}(2,3)--(2,1)--(3,1)--(3,3) to (3.1,3.3) to[out=20,in=180] node[black,above]{\scriptsize $\Gamma$} (3.5,3.7) to [out=0,in=180] (3.9,3.3) to [out=0,in=160] (4,3);
		\draw [ultra thick] (4,3.5)--(4,0)--(1,0);
		\draw [ultra thick] (0,0)--(0,4)--(2,4)--(2,1)--(3,1)--(3,3.5);
		\node at (1,1.8) {\scriptsize INK};
		\draw [line width=4] (3,-0.05)-- node[black,at end, anchor=north east]{\scriptsize HEATER} (4,-0.05);
		\draw [gray!30, dotted, line width=2.1] (0,-0.1)--(1,-0.1);
		\draw [gray!30, dotted, line width=2.1] (0,-0.2)-- node[black,below]{\scriptsize RESERVOIR} (1,-0.2);
		\draw [fill=white] (3,0) to[out=60,in=180] (3.5,0.5) to     [out=0,in=120]  (4,0);
		\node[text width=1] at (1.6,0.5) {\scriptsize VAPOR\\[-0.6ex]\scriptsize BUBBLE};
	\end{tikzpicture}
}\qquad%
\subfigure[Distorted jet.
\label{fig:t5}]{
	\begin{tikzpicture}[scale=0.8]
		\draw [gray!30,fill=gray!30] (4,3)--(4,0)--(0,0)--(0,3)--node[black,above]{\scriptsize AIR}(2,3)--(2,1)--(3,1)--(3,3) to (3.1,3.1) to[out=20,in=180]     (3.5,2.5) to [out=0,in=180] (3.9,2.9) to [out=0,in=160] (4,3);
		\draw [ultra thick] (4,3.5)--(4,0)--(1,0);
		\draw [ultra thick] (0,0)--(0,4)--(2,4)--(2,1)--(3,1)--(3,3.5);
		\node at (1,1.8) {\scriptsize INK};
		\draw [line width=4] (3,-0.05)-- node[black,at end, anchor=north east]{\scriptsize HEATER} (4,-0.05);
		\draw [gray!30, dotted, line width=2.1] (0,-0.1)--(1,-0.1);
		\draw [gray!30, dotted, line width=2.1] (0,-0.2)-- node[black,below]{\scriptsize RESERVOIR} (1,-0.2);
		\draw [dashed] (3,0) to[out=60,in=180] (3.5,0.5) to     [out=0,in=120]  (4,0);
		\draw [->] (3.5,0.7)--(3.5,0.2);
		\draw [gray!30,fill=gray!30] (3.5,4.5) to [out=180,in=90] (3.2,4) to [out=-90,in=110] (3.3,3.9) to [out=-70,in=-90] (3.4,2.8) to [out=-90,in=180] (3.5,2.55) to [out=0,in=90] (3.7,3.0) to [out=90,in=-110] (3.7,3.9) to [out=70,in=-110] (3.9,4.2) to [out=70,in=0] (3.5,4.5);
	\end{tikzpicture}
}%
\caption{Operating cycle of a drop-on-demand thermal inkjet printer.}
\label{fig:stampante}
\end{figure}

Different modes of operation are adopted by printing devices, but drop-on-demand printing methods are mostly preferred to continuous jet release, in precision applications.
This means that the ink inside the printing cartridge is subject to successive impulses (by means of thermal or piezoelectric actuators), and separate jets are thus ejected from a nozzle \mynote{according to the following steps, as displayed in \figref{fig:stampante}:
\begin{enumerate}[(a)]
 \item in the initial condition, the fluid and the thermal actuator are at rest, and the free surface $\Gamma$ has a known shape;
 \item an electrical pulse activates the heater, that induces the sudden formation of a vapor bubble which pushes the ink through the nozzle;
 \item the heater is switched off and the bubble collapses, creating a counter-pressure that makes part of the ink reverse into the nozzle: inertia induces the detachment of a jet from the rest of the ink;
 \item after the ejection, the nozzle is refilled by capillary forces and oscillations occur at the meniscus $\Gamma$;
 \item if the actuator is activated again before the initial configuration is restored, the shape of $\Gamma$ is perturbed;
 \item the perturbation of the free surface determines a poor control of the following jet, thus spoiling the quality of the printing.
\end{enumerate}
}

From \figref{fig:t4}-\ref{fig:t5}, we can notice that the dynamics inside the nozzle has a major impact on the quality of the following jet formation.
Therefore, in the present paper we look for a suitable control strategy to act on the physical oscillations naturally occurring at the free surface during the filling of the capillary pipe (cf.~\figref{fig:t3}), in order to shorten the transient between two drop ejections.

\mynote{The solution of a time-dependent optimal control problem is a difficult task, and it is even more difficult
      when moving geometries are involved.
      This is the case, for example, in shape optimization problems \cite{DelfourZolesio,SokZol,MP04}.
      In some situations, a reformulation on a fixed, reference domain can help reducing the computational effort (see, e.g., \cite{FPVshapeOpt,KV,LaumenRefDom}), but if big deformations are involved, this strategy loses its suitability.
      The same arguments hold also for free surface problems, like in the present case, where the motion of the contact line \cite{Shikhmurzaev,ManScar,WalkerOns,LaurainWalker,Flaplcl} is a further source of complexity.
}

\mynote{
    In all these cases,
the design and application of a control procedure to find the optimal control entail refined mathematical results and a high computational burden, both in terms of computational power and memory storage.
In order to overcome these difficulties, approximate optimization strategies have been developed, such as the receding horizon control - also called model predictive control - \cite{MWB,JYH01,NP97,RM93}, that can be adopted to find a suboptimal solution by solving a sequence of (inexact) optimization problems on portions of the complete timespan of interest.
Due to the low computational effort they require, these approaches are highly valuable in industrial applications, where quick - albeit suboptimal - solutions are very helpful in the design of innovative devices.
}

\mynote{
Following this perspective, in this work an instantaneous control approach (IC) is adopted, which is the simplest version of receding horizon control, initially introduced in \cite{origInstCtrl}.
Some theoretical results can be found in the literature, showing the stability and convergence properties of IC \cite{Hinze,HY97,LKC98}, and its relationship with other feedback control techniques \cite{CHK99,GM00,JSK97, Heink}.
Numerical experiments will show the effectiveness of the application of IC to the present problem, by reducing the natural oscillations of the free surface and shortening the transient that the flow experiences before achieving its equilibrium configuration.
}

The present paper is organized as follows.
\secref{sec:problem} is devoted to the description of the mathematical model of the physical phenomenon under investigation, and its stabilized ALE-FEM discretization.
In \secref{sec:ctrl}, the optimal control problem is introduced, and an instantaneous control algorithm is designed for its solution.
Numerical tests are reported in \secref{sec:instctrlResults}, in order to calibrate the control procedure and to show its effectiveness for the problem under inspection.

\section{Description and discretization of the state problem}\label{sec:problem}

\begin{figure}
 \centering
 \setlength{\figurewidth}{0.5\textwidth}
 \scriptsize
  \centering

\begin{tikzpicture}[scale=2.2, >=latex]
  \draw[draw=none, fill=gray!25] (0,0) -- node[left,pos=0.5]{$\Sigma$} (0,1.5) 
    to[out=-45,in=180] ++(0.5,-0.3) to[out=0,in=-135] node[below,pos=0]{$\Gamma$} ++(0.5,0.3) -- node[right,pos=0.5]{$\Sigma$} (1,0) -- node[above,pos=0.5]{$\Sigma_b$} cycle;
  \draw[draw=none, fill=gray!60]
    (1,1.5) arc (0:-180:0.5 and 0.1)
    to[out=-45,in=180] ++(0.5,-0.3) to[out=0,in=-135] ++(0.5,0.3);
  \draw[draw=none, fill=gray!40] (0.5,1.5) ellipse (0.5 and 0.1);

  \draw (0,0) -- (0,2);
  \draw (1,0) -- (1,2);
  \draw[yshift=2cm, dashed] (0,0) -- (0,0.15);
  \draw[yshift=2cm, dashed] (1,0) -- (1,0.15);

    \draw (0,1.5) to[out=-45,in=180] ++(0.5,-0.3) to[out=0,in=-135] ++(0.5,0.3);
    \draw (0,1.5) arc (180:360:0.5 and 0.1);
    \draw[dashed] (1,1.5) arc (0:180:0.5 and 0.1);
    \draw[dashed] (0.5,1.8) ellipse (0.5 and 0.1);


  \draw (0,1.4) arc (270:315:0.1) node[below,pos=0.5]{$\theta$};
  \draw (1,1.4) arc (270:225:0.1) node[below,pos=0.5]{$\theta$};

  \draw[-] (0.5,1.5)++(180:0.5) -- ++(135:-0.5);
  \draw[->] (0.5,1.5)++(0:0.5)++(45:-0.5) -- ++(45:1) node[right]{$\mathbf b$};

  \draw[->, line width=2pt] (0.5,0.8) -- node[right] {$\mathbf g$}++(0,-0.5);
  \draw[->, line width=2pt] (0.5,2.2) -- ++(0,0.5)
    node[above] {$\gamma \curv \normal$};
  \foreach \x in {0,60,...,300} {
      \draw[->, line width=2pt] (0.5,1.5)++(\x:0.5 and 0.1) -- ++(0,0.5);
  }
  \node[anchor=east] at (0,1.85) {$\gamma(\cos\theta-\cos\theta_s)$};
  \foreach \x in {0.1,0.3,...,0.9} {
      \draw[->, line width=2pt] (\x,0) -- ++(0,-0.2);
  }
  \node[anchor=north] at (0.5,-0.3) {$-\overline p\normal + \fbox{$\controlvar$}$};
  
  \node[anchor=east] (omega) at (-0.25,1.3) {$\Omega$};
  \draw (omega.south) to[out=-90, in=90] (0.2,0.75);  
\end{tikzpicture}
\qquad\qquad
%
%
\begin{tikzpicture}[scale=1.5]
		\filldraw[fill=gray!30,draw=white,opacity=0.6] (1,1) to [out=-135,in=20] (0,0.4) -- (0,0) -- (1,0) -- cycle;
		\draw[ultra thick] (1,0)--node[at start, right]{$\Sigma$}(1,1);
		\node[anchor=south west] at (0,0.4) {$\Gamma$};
		\draw (1,1) to [out=-135,in=20] (0,0.4);
		\draw (1,0.8) arc (-90:-135:0.2) node[below,pos=0.5] {$\theta$};
		\draw[-{Stealth}] (1,1) -- ++(45:0.4) node[at end, anchor=north ] {$\mathbf b$};
		\draw[-{Stealth}] (1,1) -- ++(90:0.4) node[anchor=south] {$\mathbf b_s$};
    \node (entrante) at (1,1) {$\boldsymbol\otimes$};
    \node[anchor=east,inner sep=5pt] at (entrante) {$\tangent_{\partial\Gamma}$};
		\end{tikzpicture}
 \caption{Geometrical settings of a free-surface problem.}
  \label{fig:ctrlDominioNS}
\end{figure}
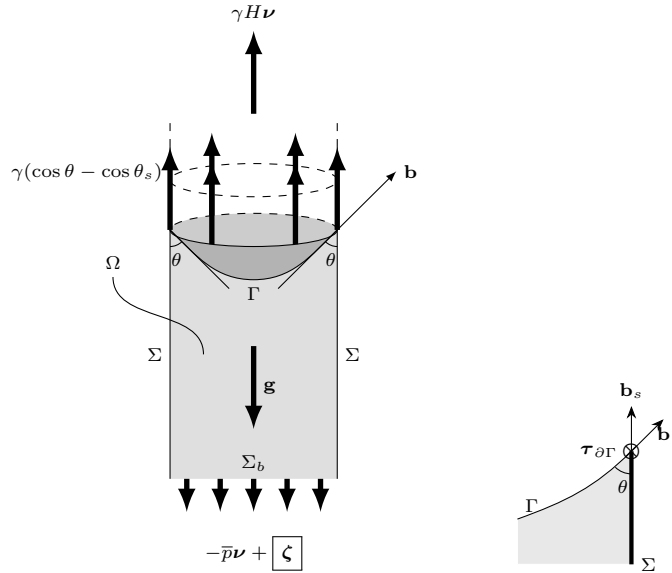

In the present work, we deal with the optimal control of a free surface flow problem.
\mynote{We consider the simple and idealized configuration depicted in \figref{fig:ctrlDominioNS}, which, anyway, contains all the main features of the addressed industrial application, such as a free surface with moving contact line and a contact force.}
An incompressible fluid lays in a cylindrical domain $\Omega\subset\mathbb R^3$, whose boundary is partitioned in a solid wall $\Sigma$ on the sides, a free surface $\Gamma$ between the fluid under investigation and a gas staying above, and a virtual, open boundary $\Sigma_b$ at the bottom, separating the region of interest from the rest of the space occupied by the fluid.
The fluid at hand is Newtonian, and we consider surface tension and capillary effects to occur at the free surface.
Therefore, a contact angle $\theta$ different from $90\degree$ can be observed at the contact line $\partial\Gamma=\overline\Gamma\cap\overline\Sigma$.

The governing equations describing the dynamics of the fluid of interest can be collected in the following time-dependent, incompressible Navier-Stokes boundary value problem on the moving domain $\Omega^t$, with boundary $\partial\Omega^t=\Gamma^t\cup\Sigma^t\cup\Sigma_b$,
\begin{equation}\label{eq:NS}
	\begin{cases}
	\partial_t\velo + (\velo\cdot \nabla)\velo - \Div\sigma = \mathbf g & \text{ in } \Omega^t, t>0,\\
	\Div\velo = 0 & \text{ in } \Omega^t, t>0,\\
	\sigma\normal\cdot\tangent = 0,\quad \sigma\normal\cdot\normal  + \gamma H= 0 & \text{ on } \Gamma^t,t>0,\\
  \velo\cdot\normal = \domvel\cdot\normal & \text{ on } \Gamma^t,t>0,\\
	\velo\cdot\normal = 0,\quad (\sigma\normal + \beta\velo+ \gamma(\cos\theta-\cos\theta_s)\delta_{\partial\Gamma}\,\mathbf b_s)\cdot\tangent = 0 & \text{ on }\Sigma^t,t>0,\\
  \sigma\normal = -\overline p\normal + \controlvar & \text { on }\Sigma_b,t>0,\\
  \velo = \velo^0 & \text{ in }\Omega^0, t=0,
	\end{cases}
\end{equation}
\mynote{where pressures $p,\overline p$ and stresses $\sigma,\controlvar$ are rescaled w.r.t. to the density $\rho$ of the fluid (e.g., $p=\widetilde p/\rho$, where $\widetilde p$ is the physical pressure).}
In system \eqref{eq:NS}, $\tangent$ is a generic vector tangent to the boundary, $\sigma=\nu\left(\nabla\velo+\nabla\velo^T\right)-p\eye$ is the (rescaled) stress tensor, $\eye$ being the identity tensor, $\mathbf g=-g\mathbf e_3$ is the gravity force, $\mathbf e_3$ being the upwards vertical vector of the canonical basis $\{\mathbf e_i\}_{i=1}^3$, $\gamma$ is the surface tension coefficient, $\curv$ is the total curvature of $\Gamma^t$, $\beta$ is the friction coefficient on $\Sigma^t$, and $-\overline p\normal +\controlvar$ is an external stress applied on $\Sigma_b$, instrumentally expressed by the sum of an hydrostatic component $-\overline p\normal$ and a generic perturbation $\controlvar$.
\mynote{The quantity $\controlvar$ will be taken as the control variable of the problem introduced in \secref{sec:ctrl}.
Therefore, in the uncontrolled case (i.e.~$\controlvar$=0) the equilibrium height of the capillary pipe will be prescribed by $\overline p$.}
The distribution $\delta_{\partial\Gamma}$ is defined as
\begin{equation}\label{eq:delta}
  \langle\delta_{\partial\Gamma}, \varphi\rangle = \int_{\partial\Gamma^t}\varphi\, d\lambda, \qquad \text{for any smooth function }\varphi,
\end{equation}
with $\lambda$ denoting the $(d-2)-$dimensional Lebesgue measure on $\partial\Gamma^t$.
\footnote{In the rest of the present paper, Lebesgue measure will be understood in all the integrals, for any dimension.}
It is worth to point out that the generalized Navier boundary condition is imposed on the solid wall $\Sigma^t$, relating the tangential stress exerted on the fluid with the discrepancy between the current contact angle $\theta$ and its equilibrium value $\theta_s$.
This condition was introduced in \cite{GNBC}, and its employment to take into account a moving contact line has been quite established in the literature of the last decade \cite{QianOns,WalkerOns,GrasselliOns,MovingCL}.

Before addressing the optimal control of problem \eqref{eq:NS}, we complete its description by discussing its numerical discretization: this is the subject of the \secref{sec:discr}.

\subsection{Discretization of the state problem}\label{sec:discr}

\mynote{In the present work, we adopt an ALE-FEM discretization of the system, inspired by \cite{Gerbeau} and discussed in \cite{MovingCL}.}
Introducing a time grid $\{t^\tn\}_{n=0}^N$ with uniform step length $\Delta t$ on the timespan $[0,T]$, a sequence $\{\Omega^\tn\}_{n=0}^N$ of domains is defined, representing the configurations of the moving domain $\Omega^t$ at the grid instants.
Correspondingly, a sequence of ALE maps of the form $\mathcal A_{n-1,n}=\eye+\Delta t\ \domvel^\tnm1$ can be introduced, mapping the domains through the different time steps: $\Omega^\tn=\mathcal A_{n-1,n}(\Omega^\tnm1)$ \cite{ALE}.
Employing the same map, we can recursively define the velocity and pressure variational spaces as follows, starting from the initial ones $\veloset^0=\{\test\in[H^1(\Omega^0)\cap H^1(\Gamma^0)]^d\st\test\cdot\normal=0\text{ on }\Sigma^0\}$ and $\pset^0=L^2(\Omega^0)$:
\begin{align}
 \veloset^\tn &= \{\test\in[H^1(\Omega^\tn)\cap H^1(\Gamma^\tn)]^d \st \test\circ\mathcal A_{n-1,n}\in\veloset^\tnm1\},\\
 \pset^\tn &= \{\ptest\in L^2(\Omega^\tn) \st \ptest\circ\mathcal A_{n-1,n}\in\pset^\tnm1\}.
\end{align}

The domain velocity $\domvel^\tn$ brings additional transport terms into the differential system \eqref{eq:NS}, whose time discretization in weak form reads as follows:
given $\velo^0$, for each $n=0,\dots,N-1$, find $(\velo^\tnp{1},p^\tnp{1})\in V^\tnp{1}\times P^\tnp{1}$ such that, $\forall(\test,\ptest)\in V^\tn\times P^\tn$,
\begin{equation}\label{eq:tdiscr}\begin{aligned}
  &\frac{1}{\Delta t}(\velo^\tnp{1},\test)_{\Omega^\tnp{1}} + a^\tnp{1}(\velo^\tnp{1},\test) + b^\tnp{1}(\test,p^\tnp{1}) \\
  &\quad- b^\tnp{1}(\velo^\tnp{1},\ptest) + c_{ALE}^\tnp{1}(\velo^\tn,\domvel^\tn,\velo^\tnp{1},\test) \\
  &\quad+ s^\tnp1(\velo^\tn,\domvel^\tn,\velo^\tnp{1},\test) + \Delta t\;S_\Gamma^\tnp1(\velo^\tnp1,\test) \\
  &\qquad= \frac{1}{\Delta t}(\velo^\tn,\test)_{\Omega^\tn} + F^\tnp{1}(\test;\controlvar),
\end{aligned}\end{equation}
where the different forms are defined as follows:
\begin{equation}\begin{split}
  (\cdot,\cdot)_{\Omega^\tnp1} &\quad\text{is the $L^2$ inner product on } \Omega^\tnp1, \\
  a^\tnp1(\velo,\test) &= \left(\frac{\nu}{2}(\nabla\velo+\nabla\velo^T),\nabla\test+\nabla\test^T\right)_{\Omega^\tnp1} + \int_{\Sigma^\tnp1} \beta\velo\cdot\test, \\
  b^\tnp1(\test,\ptest) &= -(\Div\,\test,\ptest)_{\Omega^\tnp1}, \\
  c_{ALE}^\tnp1(\other,\domvel,\velo,\test) &= \left([(\other-\domvel)\cdot\nabla]\velo,\test\right)_{\Omega^\tnp1} - \left(\Div(\domvel)\velo,\test\right)_{\Omega^\tnp1}, \\
  s^\tnp1(\other,\domvel,\velo,\test) &= \frac{1}{2}\left(\Div(\other)\velo,\test\right)_{\Omega^\tnp1} - \frac{1}{2}\int_{\Gamma^\tnp1}(\other-\domvel)\cdot\normal\ \velo\cdot\test,\\
  S_\Gamma^\tnp1(\velo,\test) &=\frac{1}{2}\int_{\Gamma^\tnp1}\gamma\normalcomp_3^2\left(\normalcomp_1\partial_3\frac{\velo\cdot\normal}{\normalcomp_3}-\normalcomp_3\partial_1\frac{\velo\cdot\normal}{\normalcomp_3}\right)\left(\normalcomp_1\partial_3\frac{\test\cdot\normal}{\normalcomp_3}-\normalcomp_3\partial_1\frac{\test\cdot\normal}{\normalcomp_3}\right),\\
  F^\tnp1(\test;\controlvar) &= (\mathbf g,\test)_{\Omega^\tnp1} + \int_{\Sigma_b}\controlvar\cdot\test 
  - \int_{\Gamma^\tnp1}\gamma\Div_\Gamma\test + \int_{\partial\Gamma^\tnp1} \gamma\test\cdot\mathbf b_s \cos\theta_s.
\end{split}\end{equation}
We point out that the ALE stabilization $s^\tnp1$ is standard and well-established in the literature \cite{Gerbeau,GerbeauMHD,Temam1977}, whereas the consistent stabilization form $S_\Gamma^\tnp1$ \mynote{has been introduced in \cite{MovingCL}.
Its role is to damp the possible spurious oscillations of the free surface in the case of an explicit treatment of the geometry, that would require to employ very short time steps:
the introduction of $S_\Gamma^\tnp1$ loosens this strong stability restriction on $\Delta t$ and, thus, significantly reduces the overall computational effort.}

Indeed, in the present work an explicit treatment of the geometry is considered, that is, $\Omega^\tnp{1}$ is defined as the image of the known, previous domain $\Omega^\tn$ through a map of the form $\mathcal A_{n,n+1}=\eye+\Delta t\,\domvel^\tn$.
The domain velocity $\domvel^\tn$ is defined as an harmonic extension of the normal fluid velocity at the free surface, namely the weak solution to the following problem:
\begin{equation}\label{eq:ALEpb}\begin{cases}
 \Delta\domvel^\tn=\mathbf 0 & \text{ in }\Omega^\tn,\\
 \domvel^\tn\cdot\normal=\velo^\tn\cdot\normal,\quad\partial_\normal\domvel\cdot\tangent=0 & \text{ on }\Gamma^\tn,\\
 \domvel^\tn\cdot\normal=0,\quad\partial_\normal\domvel\cdot\tangent=0 & \text{ on }\Sigma^\tn,\\
 \domvel^\tn=\mathbf 0 & \text{ on }\Sigma_b.
\end{cases}\end{equation}

\begin{rmrk}\label{rmrk:mappatest}
  In the terms involved in \eqref{eq:tdiscr}, the domain of integration does not always coincide with the domain of definition of the integrands: for example, $\test$ is defined on $\Omega^\tn$, but it appears in integrals over $\Omega^\tnp{1}=\mathcal A_{n,n+1}(\Omega^\tn)$.
  In order to keep a light notation, a change of variables via ALE mapping is understood in case this discordance occurs: e.g.\ $\int_{\Omega^\tnp{1}}\test$ actually means $\int_{\Omega^\tnp{1}}\test\circ\mathcal A_{n,n+1}^{-1}$.
\end{rmrk}

For the semi-discrete problem \eqref{eq:tdiscr}, we introduce also a space discretization, by defining a mesh $\mathcal T_h^{(0)}$ on the initial domain $\Omega^{(0)}$ and corresponding piecewise linear finite element spaces $\veloset_h^{(0)}$ and $\pset_h^{(0)}$.
The discrete spaces $\veloset_h^\tn,\pset_h^\tn$ for the velocity $\velo_h^\tn$ and the pressure $p_h^\tn$ are recursively defined by ALE mapping, employing the piecewise linear domain velocity $\domvel_h^\tnm1\in\veloset_h^\tnm1$.
Thus, the fully discrete scheme for the state problem reads as follows:
given a piecewise linear approximation $\velo_h^0$ of $\velo^0$, for each $n=0,\dots,N-1$, find $(\velo_h^\tnp{1},p_h^\tnp{1})\in V_h^\tnp{1}\times P_h^\tnp{1}$ such that, $\forall(\test_h,\ptest_h)\in V_h^\tn\times P_h^\tn$,
\begin{equation}\label{eq:discr}\begin{aligned}
  &\frac{1}{\Delta t}(\velo_h^\tnp{1},\test_h)_{\Omega^\tnp{1}} + a^\tnp{1}(\velo_h^\tnp{1},\test_H) + b^\tnp{1}(\test_h,p_h^\tnp{1}) \\
  &\quad- b^\tnp{1}(\velo_h^\tnp{1},\ptest_h) + c_{ALE}^\tnp{1}(\velo_h^\tn,\domvel_h^\tn,\velo_h^\tnp{1},\test_h) \\
  &\quad+ s^\tnp1(\velo_h^\tn,\domvel_h^\tn,\velo_h^\tnp{1},\test_h) + \Delta t\;S_\Gamma^\tnp1(\velo_h^\tnp1,\test_h) + s_p^\tnp1(p_h^\tnp1,\ptest_h) \\
  &\qquad= \frac{1}{\Delta t}(\velo_h^\tn,\test_h)_{\Omega^\tn} + F^\tnp{1}(\test_h;\controlvar),
\end{aligned}\end{equation}
where the additional Brezzi-Pitk\"aranta stabilization form
\begin{equation}
 s_p^\tnp1(p_h^\tnp1,\ptest_h) = C_sh^2\sum_{K\in\mathcal T_h^\tn}\int_K\nabla p_h^\tnp1\cdot\nabla\ptest_h
\end{equation}
is required by the choice of piecewise linear finite elements for both velocity and pressure \cite{stabP1P1}.

\begin{rmrk}[Friction coefficient]
 As pointed out in \cite{Yamamoto,MovingCL,GanesanTobiska,Renardybeta,Zaleskibeta}, the friction coefficient $\beta$ should be related to a {\em slip length} corresponding to the mesh size $h_3$ in the direction of the wall - in this case, vertical.
 In particular, the physical behavior of the system is correctly reproduced if a discrete coefficient $\beta_h$ is employed in \eqref{eq:discr}, defined in terms of an adimensional, mesh-independent parameter $\chi$, according to the relation
 \begin{equation}\label{eq:betah}
  \beta_h = \frac{\mu}{\chi\,h_3}.
 \end{equation}
\end{rmrk}

\section{Optimal control problem}\label{sec:ctrl}

The leading application described in the introduction inspired the study of the following optimal control problem, for the differential system \eqref{eq:NS}.
The goal is to optimally drive the evolution of a free-surface flow to a desired configuration in a finite timespan $[0,T]$, by acting on the stress at the open boundary $\Sigma_b$.
This aim can be formulated in mathematical terms as follows:
\begin{equation}
  \label{eq:optNStime}
  \begin{gathered}
  \text{Find }\controlvar^* = \argmin_{\controlvar\in L^2(\Sigma_b)}\left( \int_0^T\!\!\!\int_{\Omega^t}j_\Omega(\velo(\cdot,t))+\int_0^T\!\!\!\int_{\Gamma^t}j_\Gamma(\velo(\cdot,t))+\frac{\lambda}{2}\int_0^T\!\!\!\int_{\Sigma_b}|\controlvar|^2 \right),\\
  \text{subject to (\ref{eq:NS})},
\end{gathered}\end{equation}
where $j_\Omega$ and $j_\Gamma$ are generic objective functional densities, depending on the fluid flow, and $\lambda$ is a penalty parameter scaling the Tikhonov regularization term.

As pointed out in the introduction, the direct solution of the time-dependent control problem \eqref{eq:optNStime} can be quite a difficult and computationally expensive task.
An interesting technique aiming at reducing this computational burden is represented by the instantaneous control approach \cite{origInstCtrl,Hinze}.
The idea of this approach is to exploit the time discretization of the state problem, addressing a minimization problem for each time step.
Thus, in the present section we can consider the semi-discrete scheme \eqref{eq:tdiscr} as the state problem: the whole discussion holds without modifications if the fully discrete problem \eqref{eq:discr} is accounted for.

At each time $t^\tnp{1}$, we consider the following optimization problem: Find
\begin{gather}\label{eq:optNStdiscr}\begin{aligned}
  \controlvar^\tnp{1}&=\argmin_{\controlvar\in\controlset}J^\tnp1(\velo^\tnp1,\controlvar)\\
  &=\argmin_{\controlvar\in\controlset}\left( \int_{\Omega^\tnp{1}}j_\Omega(\velo^\tnp1) + \int_{\Gamma^\tnp1}j_\Gamma(\velo^\tnp1)+ \frac{\lambda}{2}\int_{\Sigma_b}|\controlvar|^2 \right),
  \end{aligned}\\
  \text{ subject to (\ref{eq:tdiscr})}.
\end{gather}

At this stage, it is relevant to point out that the solution to \eqref{eq:optNStdiscr} is suboptimal w.r.t.\ problem \eqref{eq:optNStime}.
Nevertheless, we will see from the numerical results of \secref{sec:instctrlResults} that it can provide effective control functions.

In order to design an optimization algorithm for problem \eqref{eq:optNStdiscr}, we need its optimality conditions, that can be derived by introducing the following Lagrangian functional, for each time step $t^\tnp1$:
\begin{equation}\label{eq:lagrangiana}\begin{aligned}
  &\mathcal L^\tnp{1} (\velo^\tnp1,p^\tnp1,\veloadj^\tn,\padj^\tn,\controlvar)\\
  &\quad= \int_{\Omega^\tnp1}j_\Omega(\velo^\tnp1)+\int_{\Gamma^\tnp1}j_\Gamma(\velo^\tnp1)+\frac{\lambda}{2}\int_{\Sigma_b}|\controlvar|^2 \\
  &\qquad+ \frac{1}{\Delta t}\int_{\Omega^\tn}\velo^\tn\cdot\veloadj^\tn+F^\tnp1(\veloadj^\tn;\controlvar)\\
  &\qquad- A^\tnp1(\velo^\tnp1,p^\tnp1,\veloadj^\tn,\padj^\tn;\velo^\tn,\domvel^\tn),
\end{aligned}\end{equation}
where the form $A^\tnp1$ collects all the terms at the left-hand side of \eqref{eq:tdiscr}, with the adjoint variables $\veloadj^\tn,\padj^\tn$ in place of the test functions $\test,\ptest$.

As usual, requiring the stationarity of the Lagrangian w.r.t.\ the adjoint variables $\veloadj^\tn,\padj^\tn$, we can retrieve the state problem \eqref{eq:tdiscr}.
On the other hand, imposing the stationarity of $\mathcal{L}$ w.r.t.\ the state variables, that is $\partial_{(\velo^\tnp1,p^\tnp1)}\mathcal L[(\adjtest,\padjtest)]=0$ for any $\adjtest\in\veloset^\tnp1$ and $\padjtest\in\pset^\tnp1$, yields the following adjoint problem:
Find $(\veloadj^\tn,\padj^\tn)\in V^\tn\times P^\tn$ such that, $\forall(\adjtest,\padjtest)\in V^\tnp1\times P^\tnp1$,
\begin{equation}\label{eq:adjNStdiscr}\begin{aligned}
  &\frac{1}{\Delta t}(\adjtest,\veloadj^\tn)_{\Omega^\tnp{1}} + a^\tnp{1}(\adjtest,\veloadj^\tn) + b^\tnp{1}(\veloadj^\tn,\padjtest) \\
  &\quad- b^\tnp{1}(\adjtest,\padj^\tn) + c_{ALE}^\tnp{1}(\velo^\tn,\domvel^\tn,\adjtest,\veloadj^\tn) \\
  &\quad+ s^\tnp1(\velo^\tn,\domvel^\tn,\adjtest,\veloadj^\tn,) + \Delta t\ S_\Gamma^\tnp1(\adjtest,\veloadj^\tn)\\
  &\qquad= \int_{\Omega^\tnp1}j_\Omega'(\velo^\tnp1)[\adjtest] + \int_{\Gamma^\tnp1}j_\Gamma'(\velo^\tnp1)[\adjtest].
\end{aligned}\end{equation}

\begin{rmrk}
  One can notice that the adjoint variables solve a problem over the new domain $\Omega^\tnp1$, but belong to the spaces $\veloset^\tn, \pset^\tn$ associated to the old domain $\Omega^\tn$.
  \mynote{This is a consequence of the Lagrangian approach based on functional \eqref{eq:lagrangiana}}.
\end{rmrk}

Neglecting for the sake of simplicity the stabilization terms, \eqref{eq:adjNStdiscr} can be rewritten in strong form as: Find $(\veloadj^\tn,\padj^\tn)\in V^\tn\times P^\tn$ such that, when mapped on the new domain $\Omega^\tnp1$, they fulfill the following system of differential equations:
\begin{equation}\label{eq:adjNStdiscrStrong}\begin{cases}
 -\Div\sigmaadj^\tn - (\velo^\tn-\domvel^\tn)\cdot\nabla\veloadj^\tn + \left(\frac{1}{\Delta t}-\Div\velo^\tn\right)\veloadj^\tn \\
 \qquad\qquad= j_\Omega'(\velo^\tnp1) & \text{ in }\Omega^\tnp1,\\
 \Div\veloadj^\tn = 0 & \text{ in }\Omega^\tnp1,\\
 \sigmaadj^\tn\normal + (\velo^\tn-\domvel^\tn)\cdot\normal\,\veloadj^\tn = j_\Gamma'(\velo^\tnp1) & \text{ on }\Gamma^\tnp1,\\
 \veloadj^\tn\cdot\normal=0 \qquad\text{ and }\\
 \left[\sigmaadj^\tn\normal + (\velo^\tn-\domvel^\tn)\cdot\normal\,\veloadj^\tn + \beta\veloadj^\tn\right]\cdot\tangent = 0 & \text{ on }\Sigma^\tnp1, \\
 \sigmaadj^\tn\normal=\mathbf 0 & \text{ on }\Sigma_b,
\end{cases}\end{equation}
where $\tangent$ is any tangent vector to the boundary and $\sigmaadj^\tn=\nu(\nabla\veloadj^\tn+(\nabla\veloadj^\tn)^T)-\padj^\tn\eye$ is the adjoint stress tensor.
From this strong form, it is evident that the adjoint problem is a {\em steady} problem, not involving any time advancement.
Indeed, since the optimization is performed separately in each time subinterval, the adjoint variables $\veloadj^\tn,\padj^\tn$ associated to each time slab $[t^\tn,t^\tnp1]$ are totally unrelated from one another.

Exploiting the Lagrangian functional $\mathcal L$ -- as in classical C\'ea's approach \cite{Cea} -- we can write
\begin{equation}
 \int_{\Sigma_b}\nabla_\controlvar J^\tnp1\cdot\delta\controlvar = \partial_\controlvar\mathcal L(\velo^\tnp1,p^\tnp1,\veloadj^\tn,\padj^\tn,\controlvar)[\delta\controlvar] = \int_{\Sigma_b}\left(\lambda\controlvar + \veloadj^\tn\right)\cdot\delta\controlvar,
\end{equation}
that is, the gradient is
\begin{equation}\label{eq:gradiente}
 \nabla_\controlvar J^\tnp1=\lambda\controlvar + \veloadj^\tn|_{\Sigma_b}.
\end{equation}

At this point, we have all the ingredients to formulate a control strategy for problem \eqref{eq:optNStime}, that is presented in \algoref{algo:optNS}.
We remark that, at each time iteration, just one step of the gradient method is performed, instead of a whole optimization loop: indeed, the underlying logic of \mynote{IC} is not to exactly solve a minimization problem at each time step, but to successively improve the approximation of the objective while marching forward in time (cf.\ \cite{Hinze}).

\begin{algorithm}\caption{Instantaneous control for problem \eqref{eq:optNStime}.}\label{algo:optNS}
\begin{algorithmic}[1]
 \STATE Given $\Omega^{(0)}\subset \mathbb R^d,\velo^{(0)}:\Omega^{(0)}\to\mathbb R^d,\controlvar^{(0)}:\Sigma_b\to\mathbb R^d$ 
 and a time step $\Delta t=T/N$:
 \FOR{$n=0$ \TO $N$}
    \STATE Solve the ALE problem \eqref{eq:ALEpb}. 
    \STATE Define $\Omega^\tnp1=(\eye+\Delta t\ \domvel^\tn)(\Omega^\tn)$.
      \STATE Solve the state problem \eqref{eq:tdiscr}. 
      \STATE Solve the adjoint problem \eqref{eq:adjNStdiscr}. 
      \STATE Update the control $\controlvar^\tnp1=\controlvar^\tn(1-\alpha\lambda)-\alpha\veloadj^\tn|_{\Sigma_b}$.\label{step:alpha}
 \ENDFOR
\end{algorithmic}
\end{algorithm}

\mynote{The computational efficiency of the IC approach can be noticed by looking at \algoref{algo:optNS}: at each time step, just one Euler iteration of the state problem and a steady adjoint problem need to be solved, thus the overall computational cost of the control procedure is comparable to that of a single solution of the uncontrolled state problem on the whole time interval.}

In the update step \ref{step:alpha} of \algoref{algo:optNS}, we can notice that a step-length parameter $\alpha$ is introduced.
\mynote{Its value will be chosen after a proper tuning, discussed in \secref{sec:alpha}}.
The outcome of the presented instantaneous control algorithm is a time-discrete, non-stationary control $\controlvar^\tn, n=0,\dots,N=T/\Delta t$, that drives the evolution of the solution towards the goal encoded in $j_\Omega, j_\Gamma$.

\section{Numerical results}\label{sec:instctrlResults}

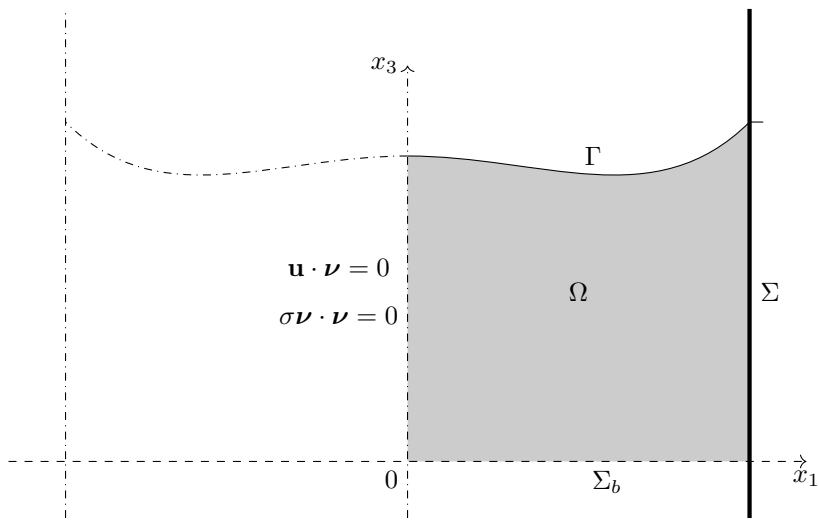
\begin{figure}
	\centering
  \begin{tikzpicture}[scale=1.5]
  \fill[black!20] (0,2.7) to[out=0,in=-135] (3,3) -- (3,0) -- (0,0);
  \draw[dashdotted,->] (0,-0.5) -- node[left, align=center]{$\velo\cdot\normal=0$\\$\sigma\normal\cdot\normal=0$} (0,3.5) node[left]{$x_3$};
  \draw[dashdotted] (-3,-0.5)--(-3,4);
  \draw[dashed,->] (-3.5,0) -- (0,0) -- node[below]{$\Sigma_b$} (3.5,0) node[below]{$x_1$};
  \draw[ultra thick] (3,4) -- (3,3) -- node[right]{$\Sigma$} (3,0) -- (3,-0.5);
  \draw (3,3) -- ++(0.12,0);
  \draw (0,2.7) to[out=0,in=-135] node[above]{$\Gamma$} (3,3);
  \draw[dashdotted] (0,2.7) to[out=180,in=-45] (-3,3);
  \node[anchor=center] at (1.5,1.5) {$\Omega$};
  \node[anchor=north east] at (0,0) {$0$};
  \end{tikzpicture}
	\caption{Axisymmetric computational domain $\Omega$ of \secref{sec:instctrlResults} (gray area).\label{fig:cylsym}}
\end{figure}

In this section, we present some numerical results obtained by the application of \algoref{algo:optNS}.
The effectiveness of the instantaneous control approach is going to be shown, and the role of the Tikhonov regularization term will be discussed.
All the numerical experiments presented here have been performed using the stabilized scheme \eqref{eq:discr}, and thus no spurious, numerical oscillation will appear at any stage of the simulation.
We point out that all the simulations are set in an axisymmetric domain $\Omega$, and therefore, the state and adjoint problem \eqref{eq:tdiscr}-\eqref{eq:adjNStdiscr} are solved in cylindrical coordinates, in half a vertical section of the domain, as depicted in \figref{fig:cylsym}.

The specific control problem addressed in this section is inspired by the ink-jet printing application described in the introduction, in which the goal is to control the natural oscillations of the free surface during the evolution of the system, and thus to shorten the transient before the attainment of the equilibrium configuration (before the ejection of the following ink jet - cf.\ \figref{fig:t3}).
To this aim, we observe that the fluid velocity represents a measure of the speed at which the system evolves from the initial to the final configuration, and that the physical oscillations of the surface are generally related to an overshooting of the equilibrium level, due to the fluid velocity reaching high values during the evolution.
Therefore, the optimal control problem under inspection can be formulated in terms of a minimization of the overall fluid velocity, and we can set the objective functions in the minimization problem \eqref{eq:optNStime} as $j_\Omega=\frac{1}{2}|\velo|^2, j_\Gamma\equiv0$.
Thus, the objective functional considered in the present section is the following:
\begin{equation}\label{eq:Jnum}
 J(\velo,\controlvar) = \frac{1}{2}\int_0^T\!\!\!\int_{\Omega^t}|\velo|^2 + \frac{\lambda}{2}\int_0^T\!\!\!\int_{\Sigma_b}|\controlvar|^2.
\end{equation}
In order to apply the instantaneous control approach, we isolate the contribution related to the time subinterval $(t^\tn,t^\tnp1)$, that reads
\begin{equation} 
 J^\tnp1(\velo^\tnp1,\controlvar^\tn) = \frac{1}{2}\int_{\Omega^\tnp1}|\velo^\tnp1|^2 + \frac{\lambda}{2}\int_{\Sigma_b}|\controlvar^\tn|^2.
\end{equation}
In all the following experiments, the control stress $\controlvar^\tn$ is chosen to be constant in space and oriented in the vertical direction, that is, it can be written in terms of a scalar control 
as $\controlvar^\tn(\mathbf x) = \control^\tn\mathbf e_3$ for any $\mathbf x\in\Sigma_b$.
With this definition, the gradient of the functional $J^\tnp1$ w.r.t.\ the scalar control $\control^\tn$ is given by (cf.~\eqref{eq:gradiente})
\begin{equation}
\nabla_{\control} J^\tnp1(\veloadj^\tn,\control^\tn) = \lambda|\Sigma_b|\control^\tn + \int_{\Sigma_b}\veloadj^\tn\cdot\mathbf e_3,
\end{equation}
where the adjoint variable $\veloadj^\tn$ is the solution of \eqref{eq:adjNStdiscr}, with $j'_\Omega(\velo^\tnp1)=\velo^\tnp1$ and $j'_\Gamma\equiv0$.
Thence, the control update step \ref{step:alpha} of \algoref{algo:optNS} reads
\begin{equation}
 \control^\tnp1 = \control^\tn(1-\alpha\lambda|\Sigma_b|) - \alpha\int_{\Sigma_b}\veloadj^\tn\cdot\mathbf e_3.
\end{equation}

To design the minimization process, we are left with two parameters to be tuned: the gradient step length $\alpha$ and the penalization coefficient 
$\lambda$.
At first, we focus on determining a suitable value for $\alpha$, and in order to decouple this tuning step from 
$\lambda$, we temporarily switch off the Tikhonov regularization term by setting $\lambda\equiv0$.
A ``rigorous'' treatment of the step length would involve a line search in the gradient direction, with the application of Armijo's rule or Wolfe's condition \cite{NocedalWright}.
However, since we perform a single minimization step for each time subinterval, the actual effectiveness of these tools would be highly limited, and thus we just consider the same single value $\alpha$ for the whole optimization procedure.

\subsection{Test case 1: vanishing contact line sources}\label{sec:tc1}

In the initial configuration of the system under inspection, a cylindrical tube is partially filled with some liquid at rest, and then capillary forces, hydrostatic pressure and gravity act together while the liquid level rises, until an equilibrium configuration is achieved.
In this first test case, we consider a static contact angle $\theta_s=90\degree$, \mynote{for which the singular contact-line term at the right-hand side of \eqref{eq:tdiscr} vanishes.}
Moreover, we recall from \eqref{eq:NS} that the following conditions are imposed at the free surface and at the bottom boundary, respectively:
\begin{equation}\label{eq:instCtrlBC}\begin{aligned}
 \sigma\normal &= \gamma\curv\normal && \text{ on }\Gamma^t,\\
 \sigma\normal &= (\overline p + \control)\mathbf e_3 && \text{ on }\Sigma_b.
\end{aligned}\end{equation}
If no control is applied (namely if $\control=0$), at the equilibrium every point of the free surface $\Gamma$ is at the same height $Z_{CL}^\infty$ of the contact line, and this height is simply prescribed by Bernoulli's theorem
\begin{equation}
 g Z_{CL}^\infty = \overline p,
\end{equation}
due to the boundary conditions \eqref{eq:instCtrlBC} and the flatness of $\Gamma$.
The parameters defining this test case are collected in \tabref{tab:setInstCtrl}.
In this regard, we point out that the value of the adimensional friction coefficient $\chi=\frac{\mu}{\beta\,h_3}$ has been chosen in order to allow the free surface to naturally oscillate around the equilibrium level for a reasoonably long time, before getting at rest.

\begin{table}
\centering
\begin{tabular}{c||c}
\begin{tabular}{c|rl}
  $\nu$ & $1.87\cdot 10^{-5}$ & m\textsuperscript{2}/s\\
  $\gamma$ & $3.91\cdot10^{-8}$ & m\textsuperscript{3}/s\textsuperscript{2}\\ 
  $\chi=\frac{\nu}{\beta\,h_3}$ & $ 5\cdot10^{-5} $ &\\
  $\theta_s$ & 90 & $\degree$ \\
  $\overline p$ & $9.81\cdot10^{-4}$ & m\textsuperscript{2}/s\textsuperscript{2}\\
\end{tabular}
&
\begin{tabular}{c|rl}
  radius & $5\cdot10^{-4}$ & m\\
  initial height & $5\cdot10^{-4}$ & m\\
  $N_1,N_3$ & $16,32$ &\\
  $\Delta t$ & $2\cdot10^{-3}$ & s\\
  $C_s$ & $0.4$ &\\
\end{tabular}
\end{tabular}
\caption{Test case 1. Physical and numerical settings.
\label{tab:setInstCtrl}}
\end{table}

\subsubsection{Tuning of the gradient step length $\alpha$}\label{sec:alpha}

\begin{figure}
 \centering
 \setlength{\figurewidth}{0.39\textwidth}
 \scriptsize
  \input{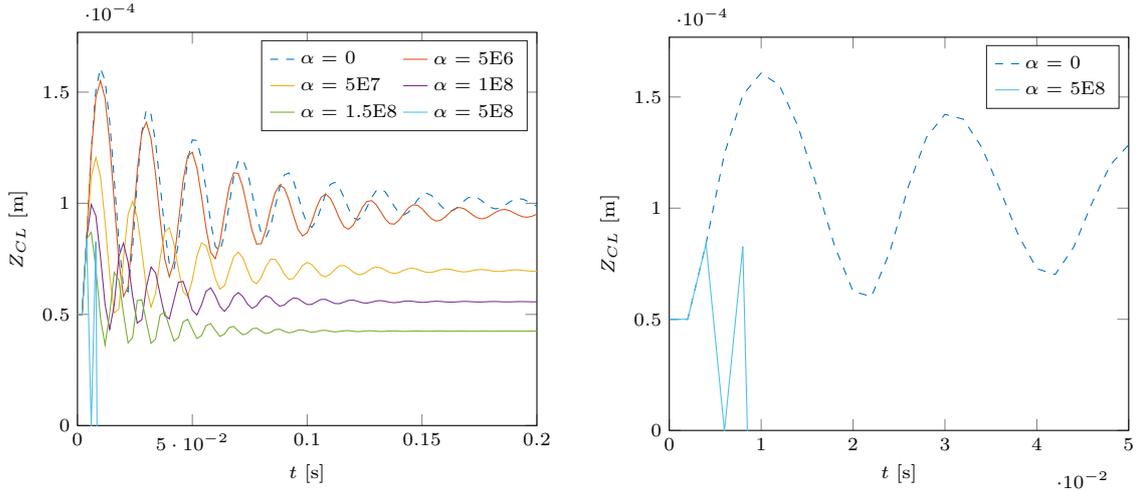}
 \quad
 \setlength{\figurewidth}{0.39\textwidth}
 \scriptsize
%
%
\definecolor{mycolor1}{rgb}{0.00000,0.44700,0.74100}%
\definecolor{mycolor2}{rgb}{0.85000,0.32500,0.09800}%
\definecolor{mycolor3}{rgb}{0.92900,0.69400,0.12500}%
\definecolor{mycolor4}{rgb}{0.49400,0.18400,0.55600}%
\definecolor{mycolor5}{rgb}{0.46600,0.67400,0.18800}%
\definecolor{mycolor6}{rgb}{0.30100,0.74500,0.93300}%
\begin{tikzpicture}

\begin{axis}[%
width=\figurewidth,
scale only axis,
xmin=0,
xmax=0.05,
ymin=0,
xlabel={$t\text{ [s]}$},
ylabel={$Z_{CL}\text{ [m]}$},
axis background/.style={fill=white},
legend style={legend cell align=left, align=left, draw=white!15!black}
]
\addplot [color=mycolor1, dashed]
  table[row sep=crcr]{%
0	5e-05\\
0.002	5e-05\\
0.004	8.40181315729763e-05\\
0.006	0.000124387838112676\\
0.008	0.000150937519073821\\
0.01	0.000160813046139019\\
0.012	0.000155162875827397\\
0.014	0.000136925946281447\\
0.016	0.000110648127776\\
0.018	8.28775302107895e-05\\
0.02	6.2584715442721e-05\\
0.022	6.02263457205663e-05\\
0.024	8.03381956314475e-05\\
0.026	0.000109541619245834\\
0.028	0.000131943119476258\\
0.03	0.000142185774433414\\
0.032	0.000140008125719443\\
0.034	0.000127467295540834\\
0.036	0.000108279359627068\\
0.038	8.77778851250394e-05\\
0.04	7.27686878806832e-05\\
0.042	7.01612133189476e-05\\
0.044	8.2357350474648e-05\\
0.046	0.000102246189337177\\
0.048	0.000119460333039501\\
0.05	0.000128526092700984\\
0.052	0.000128239881329543\\
0.054	0.000119834295805273\\
0.056	0.000106176824590935\\
0.058	9.14408531998197e-05\\
0.06	8.06479158075057e-05\\
0.062	7.83523307837515e-05\\
0.064	8.58736934617073e-05\\
0.066	9.91680062896968e-05\\
0.068	0.0001117206279464\\
0.07	0.00011906726107111\\
0.072	0.00011967508661704\\
0.074	0.000114154900834523\\
0.076	0.000104605858592684\\
0.078	9.41750611175967e-05\\
0.08	8.64928183953884e-05\\
0.082	8.459349645049e-05\\
0.084	8.92688658806372e-05\\
0.086	9.81193562188699e-05\\
0.088	0.000107030466526115\\
0.09	0.000112687694047248\\
0.092	0.000113626073950453\\
0.094	0.000110062530041917\\
0.096	0.000103469323602467\\
0.098	9.61497692607949e-05\\
0.1	9.06949423545711e-05\\
0.102	8.91602230642324e-05\\
0.104	9.20692572827788e-05\\
0.106	9.79616495347297e-05\\
0.108	0.000104202629487493\\
0.11	0.000108429252795786\\
0.112	0.000109409520894327\\
0.114	0.000107145003129024\\
0.116	0.000102634375013045\\
0.118	9.75250609961706e-05\\
0.12	9.3651630747018e-05\\
0.122	9.2428837480342e-05\\
0.124	9.42307182231687e-05\\
0.126	9.81561155044506e-05\\
0.128	0.000102495970243276\\
0.13	0.000105595804609864\\
0.132	0.000106486697554032\\
0.134	0.000105070206838361\\
0.136	0.000102006256495253\\
0.138	9.84527237795168e-05\\
0.14	9.57014338811611e-05\\
0.142	9.47385638772974e-05\\
0.144	9.5844379814813e-05\\
0.146	9.84596869497989e-05\\
0.148	0.000101465050291622\\
0.15	0.000103711571082783\\
0.152	0.000104465747992315\\
0.154	0.000103594864501802\\
0.156	0.000101525746549213\\
0.158	9.90613844754934e-05\\
0.16	9.71075505540643e-05\\
0.162	9.63576617946061e-05\\
0.164	9.70267758252221e-05\\
0.166	9.87681077876048e-05\\
0.168	0.000100843567069414\\
0.17	0.00010245852571985\\
0.172	0.000103070359297158\\
0.174	0.000102545720377223\\
0.176	0.00010115549871858\\
0.178	9.94508071317445e-05\\
0.18	9.80644100848216e-05\\
0.182	9.74864336259256e-05\\
0.184	9.78832621635072e-05\\
0.186	9.90410755043674e-05\\
0.188	0.000100471260078138\\
0.19	0.000101625418325501\\
0.192	0.000102107986757256\\
0.194	0.000101799990600051\\
0.196	0.000100870230276692\\
0.198	9.96938251813487e-05\\
0.2	9.87113535150777e-05\\
};
\addlegendentry{$\alpha\text{ = 0}$}

\addplot [color=mycolor6]
  table[row sep=crcr]{%
0	5e-05\\
0.002	5e-05\\
0.004	8.40181315729763e-05\\
0.006	-3.08716313425511e-07\\
0.008	8.27129718024399e-05\\
0.01	-0.000248037516631181\\
};
\addlegendentry{$\alpha\text{ = 5E8}$}

\end{axis}
\end{tikzpicture}%
 \caption[$Z_{CL}$ vs. time - effect of $\alpha$.]{Test case 1. Evolution of the contact line height $Z_{CL}(t)$ for different values of the discretization step: $\alpha=0$ denotes the uncontrolled case (on the right, zoom for $t\in[0,5\cdot 10^{-2}]$).}
 \label{fig:pen0}
\end{figure}

\begin{figure}
 \centering
 \setlength{\figurewidth}{0.37\textwidth}
 \scriptsize
 \input{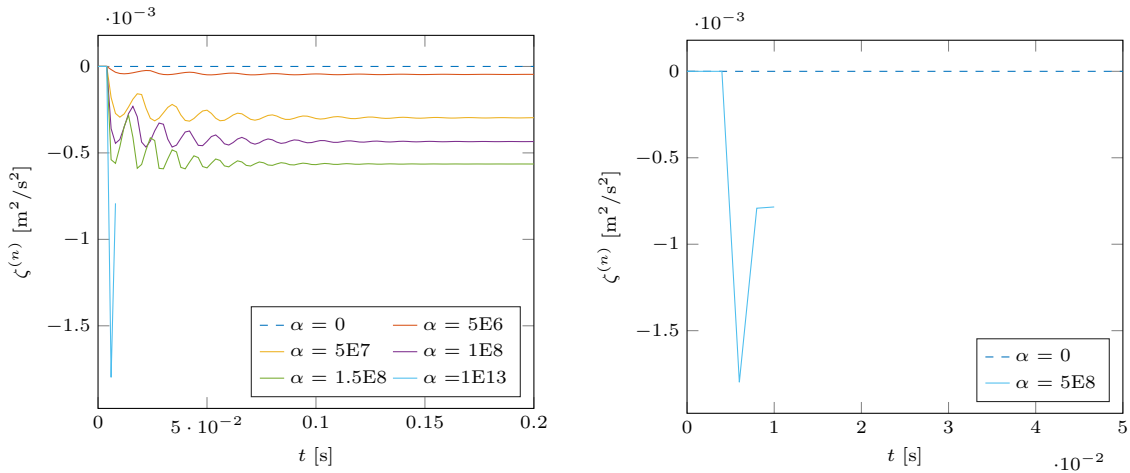}
 \quad
 \setlength{\figurewidth}{0.37\textwidth}
 \scriptsize
%
%
\definecolor{mycolor1}{rgb}{0.00000,0.44700,0.74100}%
\definecolor{mycolor2}{rgb}{0.85000,0.32500,0.09800}%
\definecolor{mycolor3}{rgb}{0.92900,0.69400,0.12500}%
\definecolor{mycolor4}{rgb}{0.49400,0.18400,0.55600}%
\definecolor{mycolor5}{rgb}{0.46600,0.67400,0.18800}%
\definecolor{mycolor6}{rgb}{0.30100,0.74500,0.93300}%
\begin{tikzpicture}

\begin{axis}[%
width=\figurewidth,
scale only axis,
xmin=0,
xmax=0.05,
xlabel={$t\text{ [s]}$},
ylabel={$\control^\tn\text{ [m\textsuperscript{2}/s\textsuperscript{2}]}$},
axis background/.style={fill=white},
legend pos=south east,
legend style={legend cell align=left, align=left, draw=white!15!black}
]
\addplot [color=mycolor1, dashed]
  table[row sep=crcr]{%
0	0\\
0.002	0\\
0.004	0\\
0.006	0\\
0.008	0\\
0.01	0\\
0.012	0\\
0.014	0\\
0.016	0\\
0.018	0\\
0.02	0\\
0.022	0\\
0.024	0\\
0.026	0\\
0.028	0\\
0.03	0\\
0.032	0\\
0.034	0\\
0.036	0\\
0.038	0\\
0.04	0\\
0.042	0\\
0.044	0\\
0.046	0\\
0.048	0\\
0.05	0\\
0.052	0\\
0.054	0\\
0.056	0\\
0.058	0\\
0.06	0\\
0.062	0\\
0.064	0\\
0.066	0\\
0.068	0\\
0.07	0\\
0.072	0\\
0.074	0\\
0.076	0\\
0.078	0\\
0.08	0\\
0.082	0\\
0.084	0\\
0.086	0\\
0.088	0\\
0.09	0\\
0.092	0\\
0.094	0\\
0.096	0\\
0.098	0\\
0.1	0\\
0.102	0\\
0.104	0\\
0.106	0\\
0.108	0\\
0.11	0\\
0.112	0\\
0.114	0\\
0.116	0\\
0.118	0\\
0.12	0\\
0.122	0\\
0.124	0\\
0.126	0\\
0.128	0\\
0.13	0\\
0.132	0\\
0.134	0\\
0.136	0\\
0.138	0\\
0.14	0\\
0.142	0\\
0.144	0\\
0.146	0\\
0.148	0\\
0.15	0\\
0.152	0\\
0.154	0\\
0.156	0\\
0.158	0\\
0.16	0\\
0.162	0\\
0.164	0\\
0.166	0\\
0.168	0\\
0.17	0\\
0.172	0\\
0.174	0\\
0.176	0\\
0.178	0\\
0.18	0\\
0.182	0\\
0.184	0\\
0.186	0\\
0.188	0\\
0.19	0\\
0.192	0\\
0.194	0\\
0.196	0\\
0.198	0\\
0.2	0\\
};
\addlegendentry{$\alpha\text{ = 0}$}

\addplot [color=mycolor6]
  table[row sep=crcr]{%
0	0\\
0.002	0\\
0.004	0\\
0.006	-0.00179799780024733\\
0.008	-0.000791685015132325\\
0.01	-0.000784483585268533\\
};
\addlegendentry{$\alpha\text{ = 5E8}$}

\end{axis}
\end{tikzpicture}%
 \caption[Control vs. time - effect of $\alpha$.]{Test case 1. Values of the control $\control^\tn$ as a function of time for different values of the discretization step $\alpha$: $\alpha=0$ denotes the uncontrolled case (on the right, zoom for $t\in[0,5\cdot10^{-2}]$).}
 \label{fig:pen0ctrl}
\end{figure}

\mynote{
Regarding the choice of the gradient step length parameter $\alpha$, it is noticed by \cite{Hinze} that a constant step is sufficient to provide relevant results, and that if a line-search-based choice of the step is performed, the computational effort is significantly increased, since two more differential problems need to be solved at each time step.
Therefore, in the present discussion we consider a constant gradient step length and we look for a suitable value for $\alpha$.
Different choices have been considered, and the resulting time evolutions of the contact line height $Z_{CL}(t)$ are shown in \figref{fig:pen0}, while the corresponding plots of the control $\control^\tn$ are reported in \figref{fig:pen0ctrl}.
}

\mynote{
We can notice that for relatively small values of $\alpha$ the controlled evolution is unsurprisingly very near to the uncontrolled case, whereas the oscillations are damped more and more, and the transient shortened, as the step length is increased.
Therefore, we can conclude that employing a larger step length is advisable, since it boosts the effect of the control and it allows to obtain a time evolution that is significantly different from that of the uncontrolled system.
However, we can see that $\alpha$ must not exceed a certain threshold: in fact, for $\alpha=5\cdot 10^8$, the control applies a negative pressure that over-contrasts the capillary rise, and the domain is emptied out -- i.e.~$Z_{CL}$ is driven to zero -- before $t=0.01$s (cf.\ \figref{fig:pen0}, right). Based on this discussion, for the case at hand, an adequate value for the step length parameter is $\alpha=1.5\cdot10^8$.
}

Comparing the different evolutions plotted in \figref{fig:pen0}-\ref{fig:pen0ctrl}, we can also notice that the final equilibrium level achieved depends on the choice of $\alpha$.
Indeed, different values of the gradient step length yield different sequences of control variables $\control^\tn$, $n=1,\dots,N$, as it can be observed in \figref{fig:pen0ctrl}.
These different histories, then, lead to different final values of the control, which are directly related to the final height of the capillary column by Bernoulli's theorem, that in presence of a nonzero final value $\control^\tN$ of the control, reads
\begin{equation}
 g Z_{CL}^\infty = \overline{p} + \control^\tN.
\end{equation}
Since, on the contrary, we want the control procedure to act only on the transient without affecting the final configuration, a properly automatic switch-off of the control would be advisable: to this aim, the Tikhonov regularization term is introduced, as discussed in the next paragraph.

\subsubsection{Role of the Tikhonov regularization term}\label{sec:lambda}

As hinted in the previous paragraph, the introduction of a nonzero penalty parameter $\lambda$ for the Tikhonov regularization term can help in preventing the control variable $\control$ from spoiling the natural equilibrium configuration of the physical system.
Indeed, looking at the definition \eqref{eq:Jnum} of the functional $J^\tnp1$, we can see that whenever an equilibrium configuration is approached, the fluid velocity $\velo^\tnp1$ gets smaller and smaller, and the penalty term involving the control becomes dominant.
As a consequence, the optimization procedure aims at switching off the control.

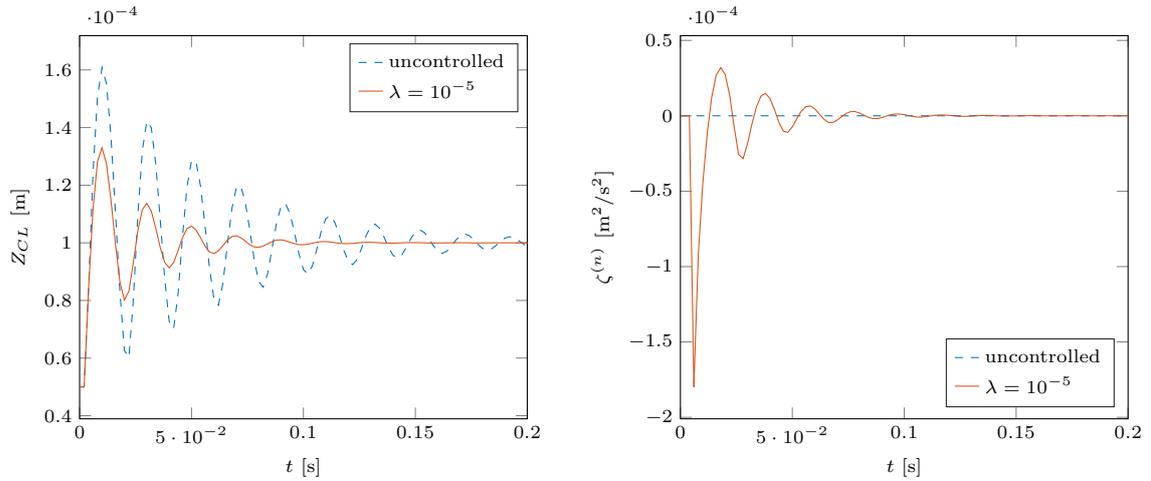
\begin{figure}
 \centering
 \setlength{\figurewidth}{0.38\textwidth}
 \scriptsize
%
%
\definecolor{mycolor1}{rgb}{0.00000,0.44700,0.74100}%
\definecolor{mycolor2}{rgb}{0.85000,0.32500,0.09800}%
\definecolor{mycolor3}{rgb}{0.92900,0.69400,0.12500}%
\definecolor{mycolor4}{rgb}{0.49400,0.18400,0.55600}%
\definecolor{mycolor5}{rgb}{0.46600,0.67400,0.18800}%
\begin{tikzpicture}

\begin{axis}[%
width=\figurewidth,
scale only axis,
xmin=0,
xmax=0.2,
xlabel={$t\text{ [s]}$},
ylabel={$Z_{CL}\text{ [m]}$},
axis background/.style={fill=white},
legend style={legend cell align=left, align=left, draw=white!15!black}
]
\addplot [color=mycolor1, dashed]
  table[row sep=crcr]{%
0	5e-05\\
0.002	5e-05\\
0.004	8.40181315729763e-05\\
0.006	0.000124387838112676\\
0.008	0.000150937519073821\\
0.01	0.000160813046139019\\
0.012	0.000155162875827397\\
0.014	0.000136925946281447\\
0.016	0.000110648127776\\
0.018	8.28775302107895e-05\\
0.02	6.2584715442721e-05\\
0.022	6.02263457205663e-05\\
0.024	8.03381956314475e-05\\
0.026	0.000109541619245834\\
0.028	0.000131943119476258\\
0.03	0.000142185774433414\\
0.032	0.000140008125719443\\
0.034	0.000127467295540834\\
0.036	0.000108279359627068\\
0.038	8.77778851250394e-05\\
0.04	7.27686878806832e-05\\
0.042	7.01612133189476e-05\\
0.044	8.2357350474648e-05\\
0.046	0.000102246189337177\\
0.048	0.000119460333039501\\
0.05	0.000128526092700984\\
0.052	0.000128239881329543\\
0.054	0.000119834295805273\\
0.056	0.000106176824590935\\
0.058	9.14408531998197e-05\\
0.06	8.06479158075057e-05\\
0.062	7.83523307837515e-05\\
0.064	8.58736934617073e-05\\
0.066	9.91680062896968e-05\\
0.068	0.0001117206279464\\
0.07	0.00011906726107111\\
0.072	0.00011967508661704\\
0.074	0.000114154900834523\\
0.076	0.000104605858592684\\
0.078	9.41750611175967e-05\\
0.08	8.64928183953884e-05\\
0.082	8.459349645049e-05\\
0.084	8.92688658806372e-05\\
0.086	9.81193562188699e-05\\
0.088	0.000107030466526115\\
0.09	0.000112687694047248\\
0.092	0.000113626073950453\\
0.094	0.000110062530041917\\
0.096	0.000103469323602467\\
0.098	9.61497692607949e-05\\
0.1	9.06949423545711e-05\\
0.102	8.91602230642324e-05\\
0.104	9.20692572827788e-05\\
0.106	9.79616495347297e-05\\
0.108	0.000104202629487493\\
0.11	0.000108429252795786\\
0.112	0.000109409520894327\\
0.114	0.000107145003129024\\
0.116	0.000102634375013045\\
0.118	9.75250609961706e-05\\
0.12	9.3651630747018e-05\\
0.122	9.2428837480342e-05\\
0.124	9.42307182231687e-05\\
0.126	9.81561155044506e-05\\
0.128	0.000102495970243276\\
0.13	0.000105595804609864\\
0.132	0.000106486697554032\\
0.134	0.000105070206838361\\
0.136	0.000102006256495253\\
0.138	9.84527237795168e-05\\
0.14	9.57014338811611e-05\\
0.142	9.47385638772974e-05\\
0.144	9.5844379814813e-05\\
0.146	9.84596869497989e-05\\
0.148	0.000101465050291622\\
0.15	0.000103711571082783\\
0.152	0.000104465747992315\\
0.154	0.000103594864501802\\
0.156	0.000101525746549213\\
0.158	9.90613844754934e-05\\
0.16	9.71075505540643e-05\\
0.162	9.63576617946061e-05\\
0.164	9.70267758252221e-05\\
0.166	9.87681077876048e-05\\
0.168	0.000100843567069414\\
0.17	0.00010245852571985\\
0.172	0.000103070359297158\\
0.174	0.000102545720377223\\
0.176	0.00010115549871858\\
0.178	9.94508071317445e-05\\
0.18	9.80644100848216e-05\\
0.182	9.74864336259256e-05\\
0.184	9.78832621635072e-05\\
0.186	9.90410755043674e-05\\
0.188	0.000100471260078138\\
0.19	0.000101625418325501\\
0.192	0.000102107986757256\\
0.194	0.000101799990600051\\
0.196	0.000100870230276692\\
0.198	9.96938251813487e-05\\
0.2	9.87113535150777e-05\\
};
\addlegendentry{uncontrolled}

\addplot [color=mycolor2]
  table[row sep=crcr]{%
0	5e-05\\
0.002	5e-05\\
0.004	8.40181315729763e-05\\
0.006	0.000111918182670065\\
0.008	0.000128370712855111\\
0.01	0.000132970826161049\\
0.012	0.000127467627929198\\
0.014	0.000114961096563189\\
0.016	9.96168615013239e-05\\
0.018	8.64005655087361e-05\\
0.02	8.01195264632593e-05\\
0.022	8.32199683810067e-05\\
0.024	9.32140351699427e-05\\
0.026	0.000104164332333842\\
0.028	0.000111536218932234\\
0.03	0.00011369864543456\\
0.032	0.000111039299789544\\
0.034	0.000105180976605511\\
0.036	9.84282285860424e-05\\
0.038	9.32081417188305e-05\\
0.04	9.13120605032524e-05\\
0.042	9.31106389154234e-05\\
0.044	9.73299161491937e-05\\
0.046	0.000101864449629931\\
0.048	0.000104952381411549\\
0.05	0.00010578069735053\\
0.052	0.000104469442932092\\
0.054	0.000101808094446189\\
0.056	9.89121118007349e-05\\
0.058	9.68491021739176e-05\\
0.06	9.62708207551258e-05\\
0.062	9.71863782789327e-05\\
0.064	9.90142267618832e-05\\
0.066	0.000100913165429903\\
0.068	0.000102171538028955\\
0.07	0.000102447679737016\\
0.072	0.000101806792830178\\
0.074	0.000100620936659856\\
0.076	9.93953787856597e-05\\
0.078	9.85800499186465e-05\\
0.08	9.8413806228203e-05\\
0.082	9.8860166169137e-05\\
0.084	9.9656756373337e-05\\
0.086	0.000100452799739679\\
0.088	0.000100957133910845\\
0.09	0.000101035100552915\\
0.092	0.000100727120606956\\
0.094	0.000100204974775428\\
0.096	9.96902763476945e-05\\
0.098	9.93692439277175e-05\\
0.1	9.93295122157151e-05\\
0.102	9.95421421188975e-05\\
0.104	9.98890947619037e-05\\
0.106	0.000100221942598092\\
0.108	0.000100421410886284\\
0.11	0.000100436504749361\\
0.112	0.00010029093733033\\
0.114	0.00010006286613475\\
0.116	9.98479045985414e-05\\
0.118	9.97224090813337e-05\\
0.12	9.97181710832077e-05\\
0.122	9.98178359502851e-05\\
0.124	9.99684877109602e-05\\
0.126	0.000100107113790294\\
0.128	0.000100185011989144\\
0.13	0.000100183545302479\\
0.132	0.000100115714645432\\
0.134	0.000100016702417424\\
0.136	9.99273648405671e-05\\
0.138	9.98788076718202e-05\\
0.14	9.9882282797125e-05\\
0.142	9.99284055875365e-05\\
0.144	9.99935603258316e-05\\
0.146	0.00010005103599836\\
0.148	0.000100081043843015\\
0.15	0.00010007704938615\\
0.152	0.000100045816701407\\
0.154	0.000100003056710616\\
0.156	9.99661072153628e-05\\
0.158	9.99475682644867e-05\\
0.16	9.99512579978251e-05\\
0.162	9.99723766733259e-05\\
0.164	0.000100000433863438\\
0.166	0.000100024149001034\\
0.168	0.000100035524562029\\
0.17	0.00010003240136615\\
0.172	0.000100018164183084\\
0.174	9.99997861760564e-05\\
0.176	9.99845797116131e-05\\
0.178	9.9977621752257e-05\\
0.18	9.99801169019945e-05\\
0.182	9.99896991039933e-05\\
0.184	0.00010000172775508\\
0.186	0.000100011463382866\\
0.188	0.000100015690377676\\
0.19	0.000100013765002848\\
0.192	0.000100007330615416\\
0.194	9.9999468455552e-05\\
0.196	9.99932425020392e-05\\
0.198	9.99906890612609e-05\\
0.2	9.99921361906717e-05\\
};
\addlegendentry{$\lambda=10^{-5}$}

\end{axis}
\end{tikzpicture}%
 \quad
 \setlength{\figurewidth}{0.38\textwidth}
 \scriptsize
%
%
\definecolor{mycolor1}{rgb}{0.00000,0.44700,0.74100}%
\definecolor{mycolor2}{rgb}{0.85000,0.32500,0.09800}%
\definecolor{mycolor3}{rgb}{0.92900,0.69400,0.12500}%
\definecolor{mycolor4}{rgb}{0.49400,0.18400,0.55600}%
\definecolor{mycolor5}{rgb}{0.46600,0.67400,0.18800}%
\begin{tikzpicture}

\begin{axis}[%
width=\figurewidth,
scale only axis,
xmin=0,
xmax=0.2,
xlabel={$t\text{ [s]}$},
ylabel={$\control^\tn\text{ [m\textsuperscript{2}/s\textsuperscript{2}]}$},
axis background/.style={fill=white},
legend pos=south east,
legend style={legend cell align=left, align=left, draw=white!15!black}
]
\addplot [color=mycolor1, dashed]
  table[row sep=crcr]{%
0	0\\
0.002	0\\
0.004	0\\
0.006	0\\
0.008	0\\
0.01	0\\
0.012	0\\
0.014	0\\
0.016	0\\
0.018	0\\
0.02	0\\
0.022	0\\
0.024	0\\
0.026	0\\
0.028	0\\
0.03	0\\
0.032	0\\
0.034	0\\
0.036	0\\
0.038	0\\
0.04	0\\
0.042	0\\
0.044	0\\
0.046	0\\
0.048	0\\
0.05	0\\
0.052	0\\
0.054	0\\
0.056	0\\
0.058	0\\
0.06	0\\
0.062	0\\
0.064	0\\
0.066	0\\
0.068	0\\
0.07	0\\
0.072	0\\
0.074	0\\
0.076	0\\
0.078	0\\
0.08	0\\
0.082	0\\
0.084	0\\
0.086	0\\
0.088	0\\
0.09	0\\
0.092	0\\
0.094	0\\
0.096	0\\
0.098	0\\
0.1	0\\
0.102	0\\
0.104	0\\
0.106	0\\
0.108	0\\
0.11	0\\
0.112	0\\
0.114	0\\
0.116	0\\
0.118	0\\
0.12	0\\
0.122	0\\
0.124	0\\
0.126	0\\
0.128	0\\
0.13	0\\
0.132	0\\
0.134	0\\
0.136	0\\
0.138	0\\
0.14	0\\
0.142	0\\
0.144	0\\
0.146	0\\
0.148	0\\
0.15	0\\
0.152	0\\
0.154	0\\
0.156	0\\
0.158	0\\
0.16	0\\
0.162	0\\
0.164	0\\
0.166	0\\
0.168	0\\
0.17	0\\
0.172	0\\
0.174	0\\
0.176	0\\
0.178	0\\
0.18	0\\
0.182	0\\
0.184	0\\
0.186	0\\
0.188	0\\
0.19	0\\
0.192	0\\
0.194	0\\
0.196	0\\
0.198	0\\
0.2	0\\
};
\addlegendentry{uncontrolled}

\addplot [color=mycolor2]
  table[row sep=crcr]{%
0	0\\
0.002	0\\
0.004	0\\
0.006	-0.000179799780024733\\
0.008	-9.18391912023659e-05\\
0.01	-4.46368998189255e-05\\
0.012	-1.1220355289822e-05\\
0.014	1.25161610571881e-05\\
0.016	2.69995790493838e-05\\
0.018	3.19607211421481e-05\\
0.02	2.73180953736294e-05\\
0.022	1.35193579997653e-05\\
0.024	-7.32719442532553e-06\\
0.026	-2.57004876252261e-05\\
0.028	-2.83665413995287e-05\\
0.03	-1.83739889472447e-05\\
0.032	-5.14630674396445e-06\\
0.034	6.09683665454486e-06\\
0.036	1.30589433066838e-05\\
0.038	1.48478489073474e-05\\
0.04	1.15423803372228e-05\\
0.042	4.30123744988513e-06\\
0.044	-4.22065281404795e-06\\
0.046	-1.01678767408128e-05\\
0.048	-1.09778033895603e-05\\
0.05	-7.38403848960223e-06\\
0.052	-1.93996179461638e-06\\
0.054	3.02452870321245e-06\\
0.056	6.05805399483996e-06\\
0.058	6.56189788399982e-06\\
0.06	4.69884243682002e-06\\
0.062	1.33519257960825e-06\\
0.064	-2.13508715192208e-06\\
0.066	-4.30711195593173e-06\\
0.068	-4.48193080754329e-06\\
0.07	-2.9545757098324e-06\\
0.072	-6.40050971635388e-07\\
0.074	1.48516191603234e-06\\
0.076	2.73083630310983e-06\\
0.078	2.81808442406612e-06\\
0.08	1.88142522870446e-06\\
0.082	3.88197257471871e-07\\
0.084	-1.03599694132851e-06\\
0.086	-1.85864309456918e-06\\
0.088	-1.8580236886749e-06\\
0.09	-1.17350827796017e-06\\
0.092	-1.77715729340538e-07\\
0.094	7.168708324262e-07\\
0.096	1.21063702867894e-06\\
0.098	1.19280798388356e-06\\
0.1	7.46473030243919e-07\\
0.102	9.54694993874225e-08\\
0.104	-4.90929773535646e-07\\
0.106	-8.04575798099646e-07\\
0.108	-7.71631360877162e-07\\
0.11	-4.6071166628877e-07\\
0.112	-3.15603330117364e-08\\
0.114	3.40966329164154e-07\\
0.116	5.32005799567333e-07\\
0.118	5.01457262531038e-07\\
0.12	2.94471002980262e-07\\
0.122	1.31769797139764e-08\\
0.124	-2.2809890810613e-07\\
0.126	-3.46873664442584e-07\\
0.128	-3.18794900529679e-07\\
0.13	-1.77646391993218e-07\\
0.132	6.8738223135978e-09\\
0.134	1.60770935703654e-07\\
0.136	2.33151731634475e-07\\
0.138	2.10613166406126e-07\\
0.14	1.161354115166e-07\\
0.142	-4.70441999443654e-09\\
0.144	-1.03671737662295e-07\\
0.146	-1.47988940045126e-07\\
0.148	-1.30064197703106e-07\\
0.15	-6.63191005861386e-08\\
0.152	1.27239786308549e-08\\
0.154	7.58739307260927e-08\\
0.156	1.02687987908023e-07\\
0.158	8.91353159963269e-08\\
0.16	4.64764521653179e-08\\
0.162	-5.18849983958713e-09\\
0.164	-4.56003462956443e-08\\
0.166	-6.17766454943616e-08\\
0.168	-5.1633670035759e-08\\
0.17	-2.30305865306684e-08\\
0.172	1.06892341297524e-08\\
0.174	3.64372575353177e-08\\
0.176	4.60982209359214e-08\\
0.178	3.86871333984887e-08\\
0.18	1.95869407321922e-08\\
0.182	-2.40558949925148e-09\\
0.184	-1.88170833032028e-08\\
0.186	-2.45411453991501e-08\\
0.188	-1.917453425555e-08\\
0.19	-6.42342151463409e-09\\
0.192	7.90045882182751e-09\\
0.194	1.83305750796922e-08\\
0.196	2.16815850375197e-08\\
0.198	1.78473355262941e-08\\
0.2	9.35573946721743e-09\\
};
\addlegendentry{$\lambda=10^{-5}$}

\end{axis}
\end{tikzpicture}%
 \caption[Effect of the penalty parameter $\lambda$ - case without contact-line force.]{Test case 1. Time evolution of the contact line height $Z_{CL}(t)$ (left) and the control variable $\control^\tn$ (right), with a nonzero penalty coefficient $\lambda$ and $\alpha=1.5\cdot10^8$.}
 \label{fig:constpen}
\end{figure}

In \figref{fig:constpen} we can see that incorporating the Tikhonov regularization into the cost functional actually yields very good results.
Indeed, the natural equilibrium level of the system is retrieved, while the time oscillations of $Z_{CL}(t)$ are effectively damped.
In addition, the duration of the transient needed to achieve the equilibrium is substantially reduced: defining $\overline t$ as the time such that
\begin{equation}\label{eq:tbar}
 |Z_{CL}(t)-Z_{CL}^\infty|< 10^{-3}\,Z_{CL}^\infty \qquad \forall t>\overline t,
\end{equation}
we have that choosing a value $\lambda=10^{-5}$ of the penalty parameter yields $\overline t=0.14$s, in opposition to a transient of more than 0.2s occurring in the uncontrolled case.
The effectiveness of this control strategy can be seen also from the plot of $\control^\tn$ w.r.t.\ time, where we see that its value quickly goes to zero, after having had a proper impact during the transient.

\vspace{1em}
\mynote{Having chosen a Tikhonov parameter $\lambda=10^{-5}$, we performed a re-calibration of the gradient step length $\alpha$.
Anyway, we found out that $\alpha=1.5\cdot10^8$ is still an effective choice.}

\subsection{Test case 2: nonzero contact-line source}\label{sec:instctrlTP}


After having discussed the application of the optimization strategy to a simple free-surface test case, we can now employ this technique to control the evolution of a system with physical parameters taken from the experimental settings reported in \cite{Yamamoto}.
The only difference with respect to the first test case, are a static contact angle $\theta_s=69.8\degree$ and the presence of a negative rescaled pressure $\overline p=-2.82\cdot10^{-2}\;\text{m}^2/\text{s}^2$ applied at the bottom boundary $\Sigma_b$.
This pressure, together with the law of capillary action, determines the following value for the final equilibrium height of the contact line:
\begin{equation}
 Z_{CL}^\infty=\frac{2\gamma\cos\theta_s}{\rho g r}-\frac{\overline p}{g} = 1.57\cdot10^{-3}\text{m}.
\end{equation}
We point out that the value of the dimensionless parameter $\chi$ is taken from \cite{MovingCL}, where the same numerical scheme as \eqref{eq:discr} is employed.

\begin{figure}
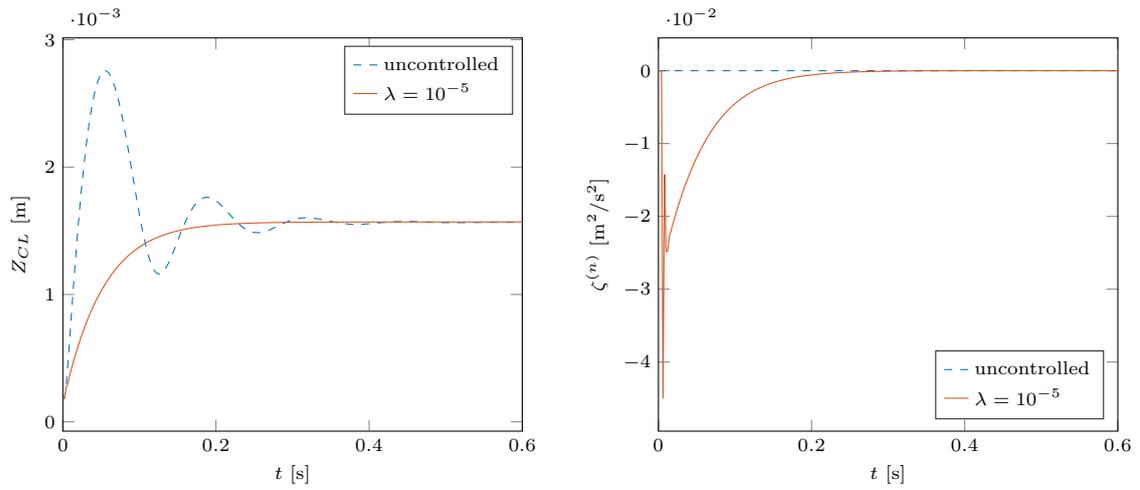

 \centering
 \setlength{\figurewidth}{0.39\textwidth}
 \scriptsize
 \input{TPconstPen}
 \quad
 \setlength{\figurewidth}{0.39\textwidth}
 \scriptsize
 \input{TPconstPen_ctrl}
 \vspace{-3mm}
 \caption[Effect of the penalty parameter $\lambda$ - case with contact-line force.]{Test case 2. Time evolution of the contact line height $Z_{CL}(t)$ (left) and the control variable $\control^\tn$ (right) in presence of a nonzero penalty parameter $\lambda$ and $\alpha=1.5\cdot10^8$.}
 \label{fig:constpenTP}
\end{figure}

\begin{figure}
 \centering
 \setlength{\fboxsep}{0pt}%
 \setlength{\fboxrule}{0.5pt}%
 \includegraphics[trim={50mm 20mm 50mm 20mm},clip,width=0.6\textwidth]{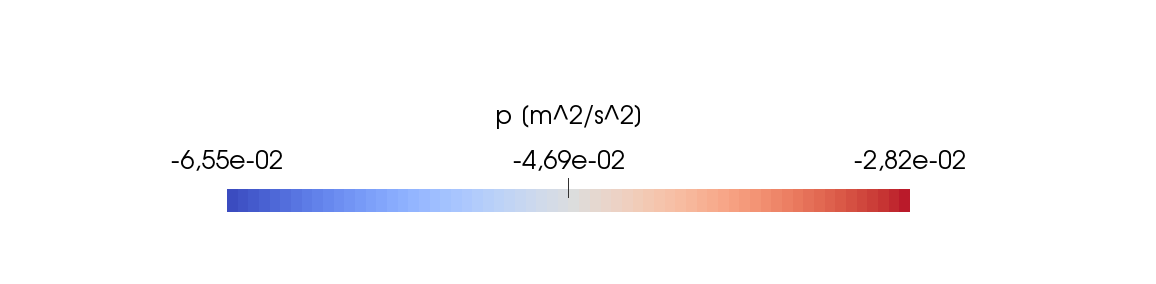}\\
 \fbox{\includegraphics[trim={10mm 1mm 126mm 0},clip,height=0.62\textwidth]{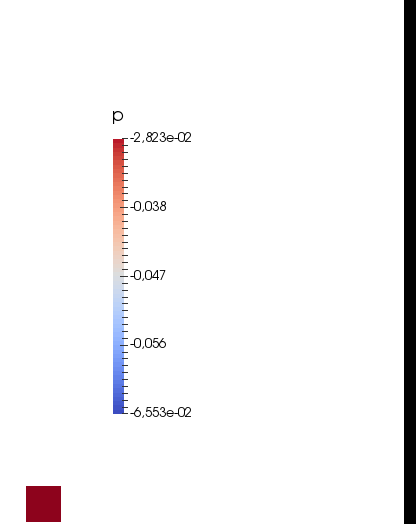}}
 \fbox{\includegraphics[trim={10mm 1mm 126mm 0},clip,height=0.62\textwidth]{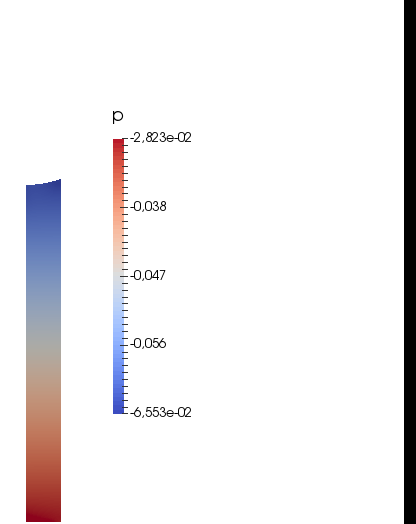}}
 \fbox{\includegraphics[trim={10mm 1mm 126mm 0},clip,height=0.62\textwidth]{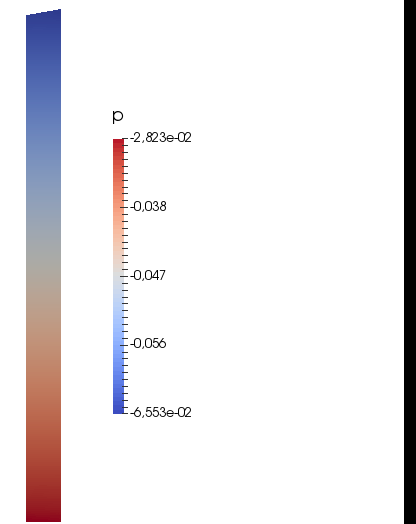}}
 \fbox{\includegraphics[trim={10mm 1mm 126mm 0},clip,height=0.62\textwidth]{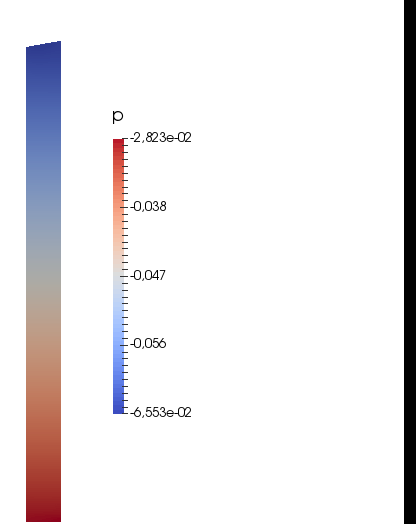}}
 \fbox{\includegraphics[trim={10mm 1mm 126mm 0},clip,height=0.62\textwidth]{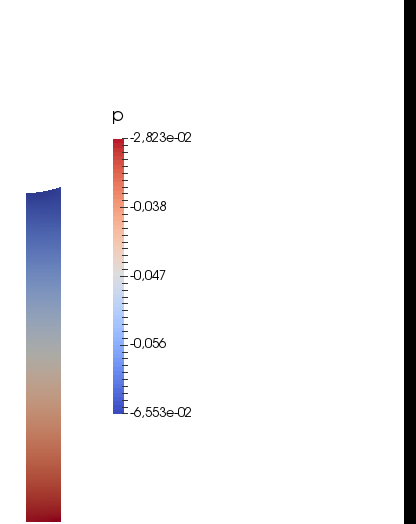}}
 \fbox{\includegraphics[trim={10mm 1mm 126mm 0},clip,height=0.62\textwidth]{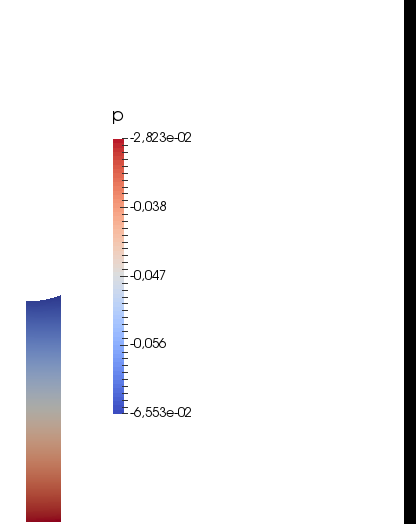}}
 \fbox{\includegraphics[trim={10mm 1mm 126mm 0},clip,height=0.62\textwidth]{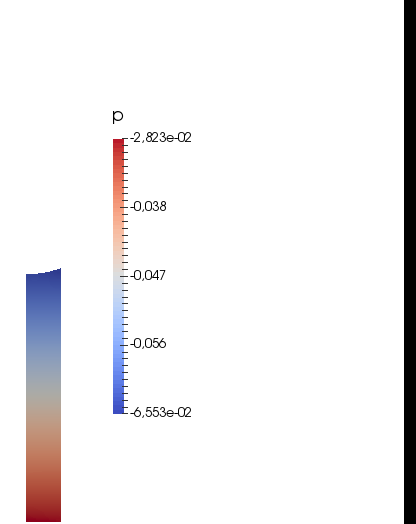}}
 \fbox{\includegraphics[trim={10mm 1mm 126mm 0},clip,height=0.62\textwidth]{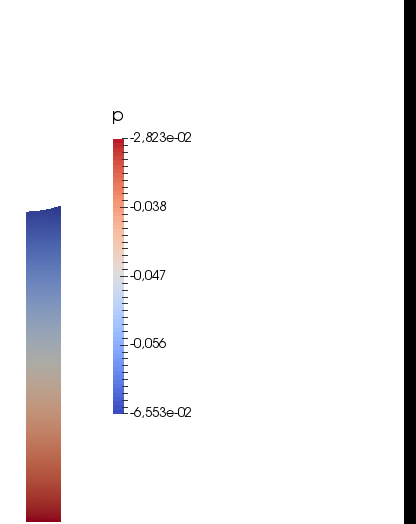}}
 \fbox{\includegraphics[trim={10mm 1mm 126mm 0},clip,height=0.62\textwidth]{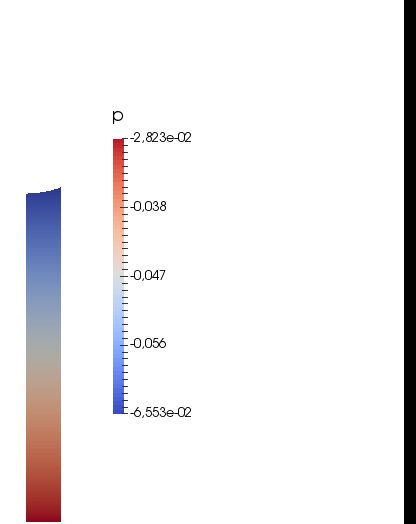}}
 \fbox{\includegraphics[trim={10mm 1mm 126mm 0},clip,height=0.62\textwidth]{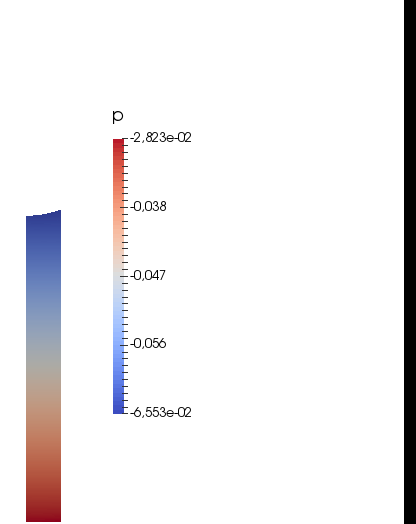}}
 \fbox{\includegraphics[trim={10mm 1mm 126mm 0},clip,height=0.62\textwidth]{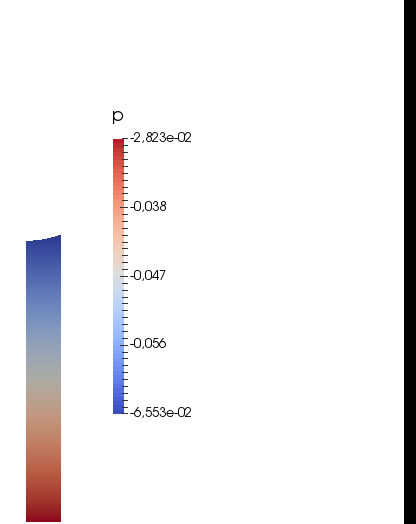}}
 \fbox{\includegraphics[trim={10mm 1mm 126mm 0},clip,height=0.62\textwidth]{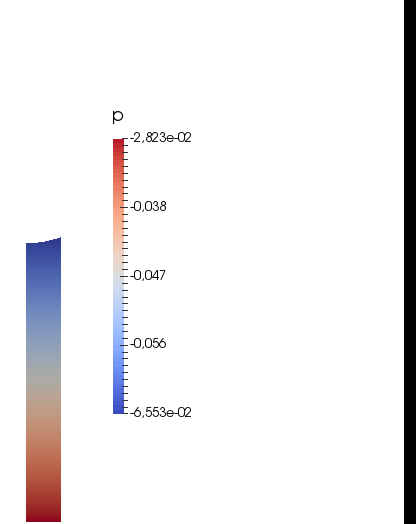}}
 \fbox{\includegraphics[trim={10mm 1mm 126mm 0},clip,height=0.62\textwidth]{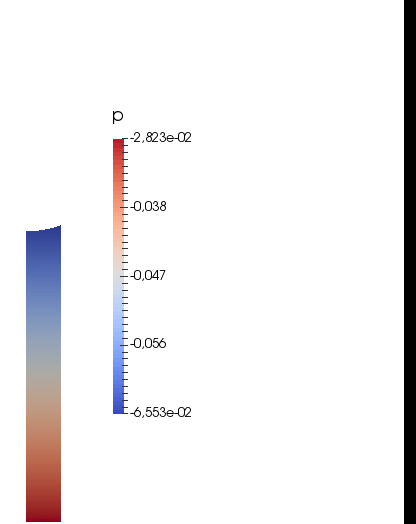}}
 \fbox{\includegraphics[trim={10mm 1mm 126mm 0},clip,height=0.62\textwidth]{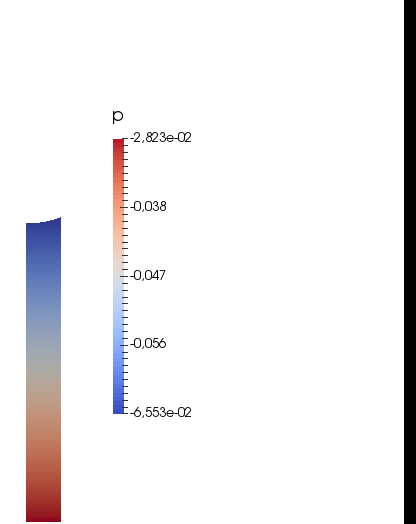}}
 \fbox{\includegraphics[trim={10mm 1mm 126mm 0},clip,height=0.62\textwidth]{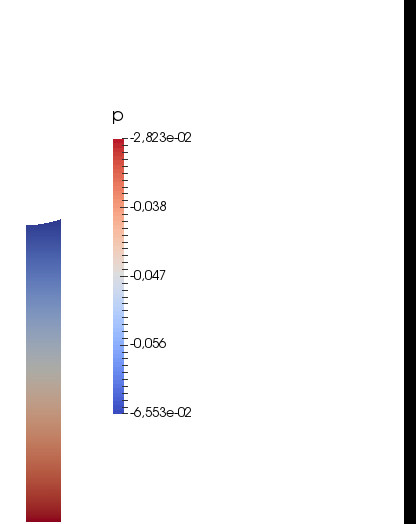}}
 \fbox{\includegraphics[trim={10mm 1mm 126mm 0},clip,height=0.62\textwidth]{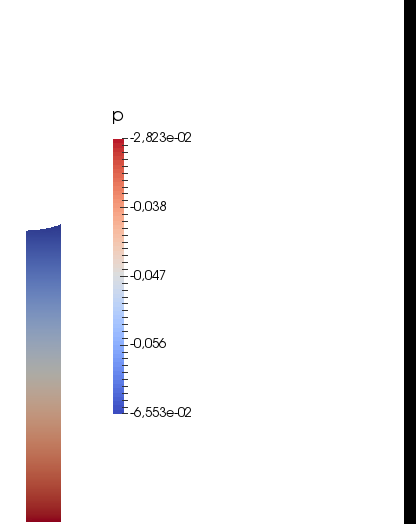}}
 \fbox{\includegraphics[trim={10mm 1mm 126mm 0},clip,height=0.62\textwidth]{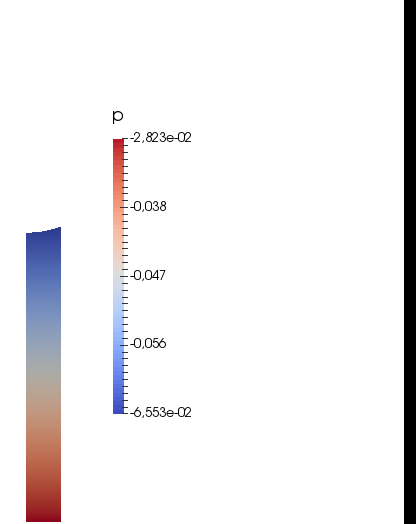}}
 \fbox{\includegraphics[trim={10mm 1mm 126mm 0},clip,height=0.62\textwidth]{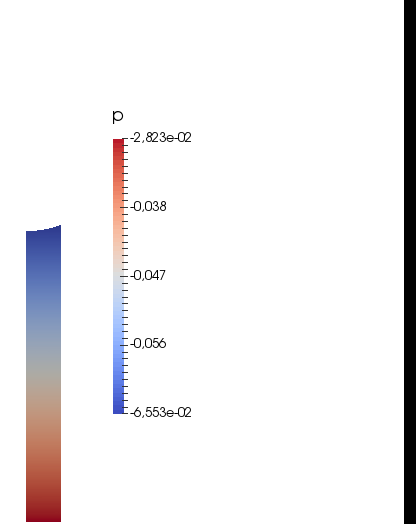}}
 \caption{Test case 2. Snapshots of the capillary rise in the uncontrolled case ($\control\equiv0$) at different time instants during the evolution with $t\in[0,0.408]$s and a sampling step of 0.024s.}
 \label{fig:snapshots}
\end{figure}

The results of this control simulation are displayed in \figref{fig:constpenTP}: since the shape of the free surface does not change significantly during the simulation (see \figref{fig:snapshots}), we can still focus on $Z_{CL}(t)$ as a measure of the overall evolution of the system.
We can see that the control strategy is extremely effective also in this case, inducing a monotone evolution towards the unperturbed natural equilibrium configuration.
Moreover, the transient is significantly shortened: the time $\overline t$ defined in \eqref{eq:tbar} to measure its duration is about $0.29$s in the controlled case, in opposition to a value $\overline t=0.54$s for the uncontrolled system.

\FloatBarrier

\section{Conclusions and future developments}\label{sec:conclusions}

In the present paper, the optimal control of the ink contained in the nozzle of an ink-jet printing device has been addressed.
The physics governing the phenomenon has been described by time-dependent incompressible Navier-Stokes equations with a free surface at the opening of the nozzle.
Surface tension has been considered at the free surface, and the dynamic contact angle imposition at the contact line has been taken into account by means of the GNBC.

\mynote{
An instantaneous control approach has been adopted to design a control strategy for the flow problem, hinging upon its time discretization.
In order to assess the resulting control procedure, the finite element method has been employed to discretize the problem and perform numerical tests inspired by the leading application.
A proper tuning of the parameters of the algorithm have allowed to provide a suitable control of the system.
Indeed, the numerical tests have shown the effectiveness and computational efficiency of the present instantaneous control approach in reducing significantly the duration of the natural transient preceding the attainment of an equilibrium configuration.
}

\mynote{
In order to enhance the treatment of the industrial problem under investigation, an interesting development on the optimal control of the free boundary could entail the actual minimization of the duration of the transient preceding the attainment of the equilibrium.
A promising perspective in this regard could be the application of the approach of \cite{Kunisch16} to the following functional:
\begin{equation}\label{eq:minKunisch}\begin{gathered}
 \controlvar = \argmin_{\controlvar\in\controlset, T\geq 0} \int_0^T\left(k+\frac{\lambda}{2}\int_{\Sigma_b}|\controlvar|^2\right)+\frac{\varepsilon}{2}\int_{\Omega^T}|\velo(\cdot,T)|^2,\\
 \text{subject to (\ref{eq:NS}) in }[0,T],
\end{gathered}\end{equation}
where $k$, $\lambda$ and $\varepsilon$ are given constants, and the last term in the functional expresses the aim of attaining equilibrium at the final time $T$.
Future work will consider this aspect, as well as the comparison of the resulting control and costs with those reported in the present paper.
In this respect, it would also be interesting to compare the IC with other model predictive control strategies, and with other shooting techniques, like the parareal method \cite{parareal} or other time-domain decomposition methods (e.g.~\cite{Heink}).
}

\section*{Acknowledgments}
The first author expresses his gratitude to Moxoff s.p.a.\ for pointing out industrial problems where the modeling and control of a moving contact line plays a crucial role, and for the financial support to his PhD activity, during which this work was started.
All the authors acknowledge the partial support of INdAM-GNCS.



\end{document}